\documentclass{amsart}

\usepackage{color}
\input xy
\xyoption{all}

\usepackage{epsfig}
\usepackage{amssymb,amsmath,amsthm,amscd}
\usepackage{fancyhdr}

\newcommand{\bC}{{\mathbb C}}
\newcommand{\bF}{{\mathbb F}}

\newcommand{\bP}{{\mathbb P}}

\newcommand{\bV}{{\mathbb V}}
\newcommand{\bW}{{\mathbb W}}
\newcommand{\bZ}{{\mathbb Z}}

\newcommand{\cK}{{\mathcal K}}
\newcommand{\cL}{{\mathcal L}}
\newcommand{\cM}{{\mathcal M}}
\newcommand{\cZ}{{\mathcal Z}}

\newcommand{\cW}{{\mathcal W}}
\newcommand{\Mbar}{\overline{\cM}}

 \DeclareMathOperator{\Ell}{Ell}
\DeclareMathOperator{\ch}{ch}

\DeclareMathOperator{\virt}{virt}

\newtheorem{theorem}{Theorem}[section]
\newtheorem{theorem/definition}{Theorem/Definition}[section]

\newtheorem{proposition}{Proposition}[section]
\newtheorem{lemma}{Lemma}[section]

\newtheorem{corollary}{Corollary}[section]

\theoremstyle{remark}

\theoremstyle{definition}

\begin{document}

\title
{On A Deformed Topological Vertex}

\author{Jian Zhou}
\address{Center of Mathatical Sciences, Zhejiang University, and
Department of Mathematical Sciences, Tsinghua University, China}
\email{jzhou@math.tsinghua.edu.cn}

\begin{abstract}
We introduce a deformed topological vertex and use it to define deformations
of the topological string partition functions of some local Calabi-Yau geometries.
We also work out some examples for which such deformations satisfy a deformed
Gopakumar-Vafa integrality and can be identified with the equivariant indices
of some naturally defined bundles on the framed moduli spaces.
\end{abstract}

\maketitle


\section{Introduction}

Since the introduction of the topological vertex in \cite{Aga-Kle-Mar-Vaf},
there have been interests in generalizing it to have more variables.
Very recently  one such generalization using Macdonald's polynomials
 has been introduced \cite{Awa-Kon}.
In this paper we will consider another one.

Recall the topological vertex involves three partitions,
when one of them is empty it is related to the leading term of
the large $N$ colored HOMFLY polynomials of the Hopf link,
and the general topological vertex can be expressed in terms of such leading terms.
We will use the whole colored HOMFLY polynomials to define the deformed vertex,
hence deformations of the partition functions of topological strings.

We define the deformed topological string partition functions simply by
replacing the topological vertices by the deformed vertices,
so their geometric meaning is missing at present.
Nevertheless, they share some important properties with the undeformed topological string partition functions
For example,
our examples indicate that the deformed topological string partition functions
satisfy a deformed Gopakumar-Vafa integrality.
We will also show that many results in \cite{ZhoCounting} and \cite{Li-Liu-Zho} can be generalized
to the deformed case,
i.e.,
we obtain some explicit expressions for some
deformed topological string partition functions which we identify with the equivariant indices
of some natural bundles on the framed moduli spaces.

The rest of the paper is arranged as follows.
In Section \ref{sec:Instanton} we compute the equivariant indices of some natural equivariant
bundles on the framed moduli spaces.
In Section \ref{sec:Schur} we collect some results on skew Schur functions
which serve as our main technical tools.
In Section \ref{sec:DTV} we will introduce the deformed topological vertex by studying the large $N$ Chern-Simons invariants of the Hopf link.
We present some product expressions for certain sums of these invariants in Section \ref{sec:Product}.
We propose in Section \ref{sec:CY} partition functions of some local Calabi-Yau
spaces using the deformed topological vertex and identify with the equivariant indices in Section \ref{sec:Instanton}.

\section{The Deformed Instanton Counting}
\label{sec:Instanton}

In this section we consider some equivariant indices of equivariant bundles
extending the ones considered in \cite{Li-Liu-Zho}.
They give us the ``deformed instonton counting functions" to be identified with
the ``deformed curved counting functions" in \S \ref{sec:CY}.

\subsection{The framed moduli spaces}
Let $M(N, k)$ denote the framed moduli space of torsion free
sheaves on $\bP^2$ with rank $N$ and $c_2=k$.
The framing means a trivialization of the sheaf restricted to the line at infinity.
In particular when $N=1$ the framed moduli spaces are the Hilbert schemes $\bC^{[k]}$.
See \cite{Nak-Yos} for details.

As proved in \cite{Nak-Yos}, $M(N,k)$ is a nonsingular variety of dimension $2Nk$.
The action of the maximal torus $T$ of $GL_N(\bC)$ together with the torus action on
$\bP^2$ induces an action on $M(N, k)$.
As shown in \cite{Nak-Yos}, the
fixed points are isolated and parameterized by $N$-tuples of
partitions $\vec{\mu}=(\mu^1, \cdots, \mu^N)$ such that $\sum_i |\mu^i|=k$.
The weight decomposition of the tangent bundle of $TM(N, k)$ at a
fixed point $\vec{\mu}$ is given by
\begin{eqnarray} \label{eqn:Weights-Tan}
&& \sum_{\alpha, \gamma=1}^N e_{\gamma}e_{\alpha}^{-1}
(\sum_{(i, j)\in \mu^\alpha} t_1^{-((\mu^{\gamma})^t_j - i)} t_2^{\mu_i^\alpha-j+1}
+\sum_{(i, j)\in \mu^\gamma} t_1^{(\mu^{\alpha})^t_j-i+1}t_2^{-(\mu_i^\gamma-j)}),
\end{eqnarray}
where $t_1, t_2\in \bC^*\times \bC^*,$ and $e_\alpha\in T$.

The space $M(N, k)$ has the following remarkable property.
Let $E$ be an equivariant coherent sheaf on it.
Even though $H^i(M(N, k), E)$ might be infinite-dimensional,
it still makes sense to define the equivariant index
and compute it by localization
(cf. \cite{Nak-Yos}):

\begin{lemma} \label{lm:Localization}
Let $E$ be an equivariant coherent sheaf on $M(N, k)$. Then
$$\chi(M(N, k), E)=\sum_{i=0}^{2Nk} (-1)^i \ch\, H^i(M(N, k), E)
=\sum_{\vec{\mu}}\ch
\left(\frac{i_{\vec{\mu}}^* E}{\wedge_{-1} T^*_{\vec{\mu}} M(N, k)}\right).$$
\end{lemma}

\subsection{Some naturally defined bundles on the framed moduli paces}

Recall $M(N, k)$ can be identified with the space of equivalent classes of tuples of linear maps
$$(B_1: V \to V; B_2: V \to V; i: W \to V; j: V \to W)$$
satisfying
$$[B_1, B_2]+ij = 0$$
and a stability condition.
Hence one gets a vector bundle $\bV_k$ over $M(N, k)$ whose fibers
are given by $V$.
This bundle is an equivariant bundle,
and its weight decomposition at a fixed point $\vec{\mu}$ is given by \cite{Nak-Yos}:
\begin{eqnarray} \label{eqn:Weights-V}
&& \bV_k = \bigoplus_{\alpha} e_{\alpha}  \sum_{(i, j) \in \mu^{\alpha}} t_1^{-i+1}t_2^{-j+1}.
\end{eqnarray}
Therefore,
the weight of ${\mbox{det}} \bV_k^*$ at the fixed point $\vec{\mu}$ is
$$\prod_{\alpha} \left( e_{\alpha}^{-|\mu^{\alpha}|}
\prod_{(i, j) \in \mu^{\alpha}} t_1^{i-1}t_2^{j-1}
\right).$$
One can has an equivariant bundle $\bW_k$ whose fibers are given by $W$.
It has the following weight decomposition:
\begin{eqnarray} \label{eqn:Weights-W}
&& \bW_k = \sum_{\alpha = 1}^N e_{\alpha}.
\end{eqnarray}
Now we take $E^m_{N, k}=K_{N, k}^{\frac{1}{2}} \otimes (\det \bV_k^*)^m$,
where $K_{N, k}$ denotes the canonical line bundle of $M(N, k)$.
By Lemma \ref{lm:Localization}, (\ref{eqn:Weights-Tan}),
(\ref{eqn:Weights-V}) and (\ref{eqn:Weights-W}) one easily gets

\begin{lemma}
We have
\begin{equation} \label{eqn:Chi}
\begin{split}
& \sum_{k \geq 0} Q^k \chi(M(N, k), \Lambda_{-e^{-t}}(\bV_k \otimes \bW_k^*) \otimes
\Lambda_{-e^{-t}}(\bV_k^* \otimes \bW_k) \otimes E_{N, k}^m)(e_1, \dots, e_N, t_1, t_2) \\
= & \sum_{\mu^{1, \dots, N}} Q^{\sum_{i=1}^N |\mu^i|}
\prod_{\alpha =1}^N
\left( e_{\alpha}^{-|\mu^{\alpha}|} \prod_{(i, j) \in \mu^{\alpha}} t_1^{i-1}t_2^{j-1} \right)^m \\
& \cdot \prod_{\alpha, \gamma=1}^N  \prod_{(i, j) \in \mu^\alpha}
\frac{1 - e^{-t} e_{\alpha}e_{\gamma}^{-1} t_1^{-i+1} t_2^{-j+1}}
{(e_{\alpha}^{-1}e_{\gamma} t_1^{-((\mu^{\gamma})^t_j - i)}
t_2^{\mu_i^\alpha-j+1})^{\frac{1}{2}}
- (e_{\alpha}^{-1}e_{\gamma}t_1^{-((\mu^{\gamma})^t_j - i)} t_2^{\mu_i^\alpha-j+1})^{-\frac{1}{2}}} \\
& \cdot \prod_{(i, j)\in \mu^\gamma} \frac{1 - e^{-t}
(e_{\gamma}e_{\alpha}^{-1} t_1^{-i+1} t_2^{-j+1})^{-1}}
{(e_{\alpha}^{-1}e_{\gamma}
t_1^{(\mu^{\alpha})^t_j-i+1}t_2^{-(\mu_i^\gamma-j)})^{\frac{1}{2}}
 - (e_{\alpha}^{-1}e_{\gamma} t_1^{(\mu^{\alpha})^t_j-i+1}t_2^{-(\mu_i^\gamma-j)})^{-\frac{1}{2}}}.
\end{split} \end{equation}
\end{lemma}

\subsection{Rank $1$ case and generalizations}

When $N=1$, $t_1= q$ an $t_2 = q^{-1}$,
we have
\begin{equation} \label{eqn:ChiVV-}
\begin{split}
& \sum_{k \geq 0} Q^k \chi((\bC^2)^{[k]}, \Lambda_{-e^{-t}}(\bV_k \otimes \bW_k^*) \otimes
\Lambda_{-e^{-t}}(\bV_k^* \otimes \bW_k) \otimes K_{1, k})(q, q^{-1}) \\
= & \sum_{\mu} Q^{|\mu|} \prod_{x \in \mu}
\frac{1 - e^{-t} q^{c(x)}}{q^{-h(x)/2} - q^{h(x)/2}} \cdot \frac{1 - e^{-t} q^{-c(x)}}
{q^{h(x)/2} - q^{-h(x)/2}}.
\end{split} \end{equation}
Since the canonical line bundle $K_{1, k}$ is equivariantly trivial,
this is equal to:
\begin{equation}
\begin{split}
& \sum_{k \geq 0} Q^k \chi(M(1, k), \Lambda_{-e^{-t}}(\bV_k \otimes \bW_k^*) \otimes
\Lambda_{-e^{-t}}(\bV_k^* \otimes \bW_k))(q, q^{-1}) \\
= & \sum_{\mu} Q^{|\mu|} \prod_{x \in \mu}
(\frac{1 - e^{-t} q^{c(x)}}{1 - q^{h(x)}} \cdot
\frac{1 - e^{-t} q^{-c(x)}}{1 - q^{-h(x)}}).
\end{split} \end{equation}

This can be generalized as follows. Let
\begin{eqnarray*}
&& W_{e^{-t}, y}(\bV_n) = \Lambda_{-e^{-t}}(\bV_n) \otimes \Lambda_{-e^{-t}}(\bV_n^*)
\otimes S_{e^{-t}y}(\bV_n) \otimes S_{e^{-t}y}(\bV_n^*).
\end{eqnarray*}
Now by (\ref{eqn:Weights-Tan}) and (\ref{eqn:Weights-V}),
\begin{equation} \label{eqn:ChiW}
\begin{split}
& \chi((\bC^2)^{[n]}, \Lambda_{-y}(T^*(\bC^2)^{[n]}) \otimes W_{e^{-t}, y}(\bV_n))(t_1, t_2) \\
= & \sum_{|\mu|=n} \prod_{(i,j)\in\mu}
\frac{(1 - e^{-t}t_1^{-i+1}t_2^{-j+1})(1 - e^{-t}t_1^{i-1}t_2^{j-1})}
{(1- t_1^{\mu^t_j-i}t_2^{-(\mu_i-j+1)})(1 - t_1^{-(\mu^t_j - i + 1)} t_2^{\mu_i -j})} \\
& \cdot  \prod_{(i,j)\in\mu}
\frac{(1- yt_1^{\mu^t_j-i}t_2^{-(\mu_i-j+1)})(1 - yt_1^{-(\mu^t_j - i + 1)} t_2^{\mu_i -j})}
{(1 - ye^{-t}t_1^{-i+1}t_2^{-j+1})(1 - ye^{-t}t_1^{i-1}t_2^{j-1})}.
\end{split}
\end{equation}
This generalizes the equivariant $\chi_y$ genera of the Hilbert schemes
as one can see by taking $e^{-t} = 0$.
It is easy to see that
\begin{equation} \label{eqn:ChiWq}
\begin{split}
& \chi((\bC^2)^{[n]}, \Lambda_{-y}(T^*(\bC^2)^{[n]}) \otimes W_{e^{-t}, y}(\bV_n))(q, q^{-1}) \\
= & \sum_{|\mu|=n} \prod_{x\in\mu}
(\frac{(1 - e^{-t}q^{c(x)})(1 - e^{-t}q^{-c(x)})}{(1 - ye^{-t}q^{c(x)})(1 - ye^{-t}q^{-c(x)})}
\frac{(1- yq^{h(x)})(1 - yq^{-h(x)})}{(1- q^{h(x))})(1 - q^{-h(x)})}).
\end{split}
\end{equation}
We will generalize this further to an analogue of the elliptic genus as follows.
Define
\begin{eqnarray*}
\Ell(T(\bC^2)^{[n]}, \bV_n)(y, p, e^{-t})
& = & \bigotimes_{n =1}^{\infty}
(\Lambda_{-yp^{n-1}}(T^*(\bC^2)^{[n]}) \otimes \Lambda_{-y^{-1}p^n}(T(\bC^2)^{[n]}) \\
&& \otimes S_{p^{n}}(T^*(\bC^2)^{[n]}) \otimes S_{p^n}(T(\bC^2)^{[n]}) \\
&& \otimes \Lambda_{-e^{-t}p^{n-1}}(\bV_n \oplus \bV_n^*)
\otimes \Lambda_{-e^{-t}p^{n}}(\bV_n \oplus \bV_n^*)^* \\
&& \otimes S_{e^{-t}yp^{n-1}}(\bV_n \oplus \bV_n^*)
\otimes S_{e^{-t}y^{-1}p^{n}}(\bV_n \oplus \bV_n^*)^*.
\end{eqnarray*}
This suggests that we should consider the supermanifolds obtained from $(\bC^2)^{[n]}$
whose odd part of the tangent spaces are $\bV_n \oplus \bV_n^*$.
By taking $p=0$,
one gets:
\begin{eqnarray*}
&& \Ell(T(\bC^2)^{[n]}, \bV_n)(y, 0, e^{-t}) \\
& = & \Lambda_{-y}(T^*(\bC^2)^{[n]}) \otimes \Lambda_{-e^{-t}}(\bV_n) \otimes \Lambda_{-e^{-t}}(\bV_n^*)
\otimes S_{e^{-t}y}(\bV_n) \otimes S_{e^{-t}y}(\bV_n^*) \\
& = & \Lambda_{-y}(T^*(\bC^2)^{[n]}) \otimes W_{e^{-t}, y}(\bV_n).
\end{eqnarray*}
by (\ref{eqn:Weights-Tan}) and (\ref{eqn:Weights-V}),
\begin{equation} \label{eqn:ChiEll}
\begin{split}
& \chi((\bC^2)^{[n]}, \Ell(T(\bC^2)^{[n]}, \bV_n)(y, p, e^{-t}))(t_1, t_2) \\
= & \sum_{|\mu|=n} \prod_{(i,j)\in\mu} \prod_{k=1}^{\infty}
\left(
\frac{(1- yp^{k-1}t_1^{\mu^t_j-i}t_2^{-(\mu_i-j+1)})(1 - yp^{k-1}t_1^{-(\mu^t_j - i + 1)} t_2^{\mu_i -j})}
{(1- p^{k-1} t_1^{\mu^t_j-i}t_2^{-(\mu_i-j+1)})(1 - p^{k-1}t_1^{-(\mu^t_j - i + 1)} t_2^{\mu_i -j})}
\right. \\
& \cdot \frac{(1- yp^kt_1^{-(\mu^t_j-i)}t_2^{\mu_i-j+1})(1 - yp^kt_1^{\mu^t_j - i + 1} t_2^{-(\mu_i -j)})}
{(1- p^k t_1^{-(\mu^t_j-i)}t_2^{\mu_i-j+1})(1 - p^k t_1^{\mu^t_j - i + 1} t_2^{-(\mu_i -j)})} \\
& \cdot \frac{(1 - e^{-t}p^{k-1}t_1^{-i+1}t_2^{-j+1})(1 - e^{-t}p^{k-1}t_1^{i-1}t_2^{j-1})}
{(1 - ye^{-t}p^{k-1}t_1^{-i+1}t_2^{-j+1})(1 - ye^{-t}p^{k-1}t_1^{i-1}t_2^{j-1})} \\
& \left. \cdot \frac{(1 - e^{-t}p^kt_1^{-i+1}t_2^{-j+1})(1 - e^{-t}p^{k-1}t_1^{i-1}t_2^{j-1})}
{(1 - ye^{-t}p^kt_1^{-i+1}t_2^{-j+1})(1 - ye^{-t}p^kt_1^{i-1}t_2^{j-1})} \right).
\end{split}
\end{equation}
It follows easily that
\begin{equation} \label{eqn:ChiEllq}
\begin{split}
& \chi((\bC^2)^{[n]}, \Ell(T(\bC^2)^{[n]}, \bV_n)(y, p, e^{-t}))(q, q^{-1}) \\
= & \sum_{|\mu|=n} \prod_{(i, j)\in \mu}
\prod_{k=1}^\infty
\left( \frac{(1-p^{k-1}yq^{h(i, j)})(1-p^{k-1}yq^{-h(i, j)})}
{(1-p^{k-1}q^{h(i,j)})(1-p^{k-1}q^{-h(i, j)})} \right. \\
& \cdot \frac{(1-p^ky^{-1} q^{h(i, j)})(1-p^ky^{-1}q^{-h(i, j)})}
{(1-p^kq^{h(i,j)})(1-p^kq^{-h(i, j)})} \\
& \cdot \frac{(1-e^{-t} p^{k-1}q^{c(i,j)})(1-e^{-t}p^{k-1}q^{-c(i, j)})}
{(1-e^{-t} p^{k-1}yq^{c(i, j)})(1-e^{-t}p^{k-1}yq^{-c(i, j)})} \\
& \left. \cdot \frac{(1-e^{-t}p^kq^{c(i,j)})(1-e^{-t}p^kq^{-c(i, j)})}
{(1-e^{-t}p^ky^{-1} q^{c(i, j)})(1-e^{-t}p^ky^{-1}q^{-c(i, j)})} \right).
\end{split}
\end{equation}
It is natural to expect that a deformed version of the equivariant DMVV conjecture holds
for such deformed elliptic genera.

\subsection{Rank $>1$ case}
When $t_1 = q$, $t_2 = q^{-1}$,
after a tedious elementary calculation one can find:
\begin{equation} \label{eqn:Rank>1}
\begin{split}
& \sum_{k \geq 0} Q^k \chi(M(N, k), \Lambda_{-e^{-t}}(\bV_k \otimes \bW_k^*) \otimes
\Lambda_{-e^{-t}}(\bV_k^* \otimes \bW_k) \otimes E_{N, k}^m)(e_1, \dots, e_N, q, q^{-1}) \\
= & \sum_{\mu^{1, \dots, N}} Q^{\sum_{i=1}^N |\mu^i|} e^{-Nt \sum_{\alpha=1}^N |\mu^{\alpha}|}
\prod_{\alpha =1}^N
\left( e_{\alpha}^{-|\mu^{\alpha}|} q^{-\kappa_{\mu^{\alpha}}/2} \right)^m \\
& \cdot \prod_{\alpha=1}^N \prod_{(i, j) \in \mu^\alpha}
\frac{e^{t/2} q^{(j-i)/2} - e^{-t/2} q^{-(j-i)/2}}{q^{h(i,j)/2} - q^{-h(i,j)/2}} \\
& \cdot \prod_{1 \leq \alpha < \gamma \leq N}  (\prod_{(i, j) \in \mu^\alpha}
\frac{1 - e^{-t} e_{\alpha}e_{\gamma}^{-1} q^{j-i}}
{1 - e_{\alpha}e_{\gamma}^{-1} q^{\mu_i^\alpha+(\mu^{\gamma})^t_j - i-j+1}} \\
& \cdot \prod_{(i, j)\in \mu^\gamma} \frac{1 - e^{-t} e_{\alpha}e_{\gamma}^{-1} q^{-(j-i)}}
{1 - e_{\alpha}e_{\gamma}^{-1} q^{-(\mu_i^\gamma+(\mu^{\alpha})^t_j-i-j+1)}}) \\
& \times \prod_{\alpha=1}^N \prod_{(i, j) \in \mu^\alpha}
\frac{e^{-t/2} q^{(j-i)/2} - e^{t/2} q^{-(j-i)/2}}{q^{h(i,j)/2} - q^{-h(i,j)/2}} \\
& \cdot \prod_{1 \leq \alpha < \gamma \leq N}  (
\prod_{(i, j)\in \mu^\alpha} \frac{1 - e^{t} e_{\alpha}e_{\gamma}^{-1} q^{j-i}}
{1 - e_{\alpha}e_{\gamma}^{-1} q^{\mu_i^\alpha+(\mu^{\gamma})^t_j-i-j+1}} \\
& \cdot \prod_{(i, j) \in \mu^\gamma}
\frac{1 - e^{t} e_{\alpha}e_{\gamma}^{-1} q^{-(j-i)}}
{1 - e_{\alpha}e_{\gamma}^{-1} q^{-(\mu_i^\gamma+(\mu^{\alpha})^t_j - i-j+1)}}).
\end{split}
\end{equation}

\section{Preliminary Results on Skew Schur Functions}
\label{sec:Schur}

In this section we collect some results on skew Schur functions.
They will be our main technical tools for the rest of the paper.

\subsection{Summation formulas for Skew Schur functions}

Recall the following identity:
\begin{eqnarray} \label{eqn:ProdSum}
&& \prod_{i, j =1}^{\infty} (1 - Qx_iy_j)^{-1}
= \exp \sum_{n=1}^{\infty} \frac{Q^n}{n}p_n(x)p_n(y),
\end{eqnarray}
where $x=(x_1, x_2, \dots)$, $y=(y_1, y_2, \dots)$,
and
$$p_n(x) = x_1^n + x_2^n + \cdots.$$
For some specializations,
the left-hand side of (\ref{eqn:ProdSum}) may not make sense,
but the right-hand side makes sense.
For example,
when $x = q^{\rho} = (q^{-\frac{1}{2}}, q^{-\frac{3}{2}}, \dots)$, $y = q^{-\rho}
= (q^{\frac{1}{2}}, q^{\frac{3}{2}}, \dots)$,
the left-hand side becomes
$$ \prod_{i, j = 1}^{\infty} (1 - Qq^{j-i}),$$
which does not make sense, while on the right-hand,
\begin{eqnarray} \label{eqn:PnQrho}
&& p_n(q^{\rho})
= q^{-\frac{n}{2}} + q^{-\frac{3n}{2}} + \cdots =
\frac{q^{-\frac{n}{2}}}{1 - q^{-n}} = - \frac{q^{\frac{n}{2}}}{1 -
q^n} = - p_n(q^{-\rho}),
\end{eqnarray}
hence the right-hand side of
(\ref{eqn:ProdSum}) is
$$\exp \left(-\sum_{n=1}^{\infty} \frac{(qQ)^n}{(q^n-1)^2}\right).$$

We will use the following identity for skew Schur functions (cf. \cite[p. 93, (1)]{Mac}):
\begin{eqnarray}
&& \sum_{\eta} s_{\eta/\mu}(x)s_{\eta/\nu}(y)
=  \exp \sum_{n = 1}^{\infty} \frac{1}{n} p_n(x)p_n(y)
\cdot \sum_{\tau} s_{\mu/\tau}(y)s_{\nu/\tau}(x),  \label{eqn:SchurSum} \\
&& \sum_{\eta} s_{\eta/\mu^t}(x)s_{\eta^t/\nu}(y)
=  \exp \sum_{n = 1}^{\infty} \frac{(-1)^{n-1}}{n} p_n(x)p_n(y)
\cdot \sum_{\tau} s_{\mu/\tau^t}(y)s_{\nu^t/\tau}(x),  \label{eqn:SchurSumT}
\end{eqnarray}
where $x=(x_1, x_2, \dots)$, $y=(y_1, y_2, \dots)$.
In particular,
\begin{eqnarray}
&& \sum_{\eta} s_{\eta}(x)s_{\eta}(y)
= \exp \sum_{n = 1}^{\infty} \frac{1}{n} p_n(x)p_n(y), \label{eqn:SchurSum1} \\
&& \sum_{\eta} s_{\eta/\mu}(x)s_{\eta}(y)
=  \exp \sum_{n = 1}^{\infty} \frac{1}{n} p_n(x)p_n(y) \cdot  s_{\mu}(y),  \label{eqn:SchurSum2} \\
&& \sum_{\eta} s_{\eta}(x)s_{\eta^t}(y)
= \exp \sum_{n = 1}^{\infty} \frac{(-1)^{n-1}}{n} p_n(x)p_n(y), \label{eqn:SchurSum1T} \\
&& \sum_{\eta} s_{\eta/\mu^t}(x)s_{\eta}(y)
=  \exp \sum_{n = 1}^{\infty} \frac{(-1)^{n-1}}{n} p_n(x)p_n(y) \cdot  s_{\mu}(y).  \label{eqn:SchurSum2T}
\end{eqnarray}

For $x = (x_1, x_2, \dots)$, $y = (y_1, y_2, \dots)$,
let $(x, y) = (x_1, y_1, x_2, y_2, \dots)$.
Then we have
\begin{eqnarray} \label{eqn:SchurSumXY}
&& \sum_{\xi} s_{\mu/\xi}(x) s_{\xi/\nu}(y) = s_{\mu/\nu}(x, y).
\end{eqnarray}

We will need the following:

\begin{lemma} \cite{ZhoCounting}
The following identity holds:
\begin{equation} \label{eqn:SchurSum3}
\begin{split}
& \sum_{\nu^1, \dots, \nu^{N}}\sum_{\eta^1, \dots, \eta^{N-1}}
\prod_{k=1}^{N} s_{\nu^k/\eta^{k-1}}(x^k)Q_k^{|\nu^k|} s_{\nu^k/\eta^{k}}(y^k) \\
= & \prod_{1 \leq k < l \leq N+1} \exp \sum_{n = 1}^{\infty}
\frac{(Q_kQ_{k+1} \cdots Q_{l-1})^n}{n}  p_n(x^k)p_n(y^{l-1}),
\end{split} \end{equation}
where $\eta^0=\eta^N = (0)$,
$x=(x^k_1, x^k_2, \dots)$, $y^k=(y^k_1, y^k_2, \dots)$.
\end{lemma}

Formulas of the type in the following Lemma appeared in \cite{Hol-Iqb-Vaf} without proof.

\begin{lemma}
For $|uv| < 1$ we have
\begin{eqnarray}
&& \sum_{\mu, \nu} u^{|\mu|} v^{|\nu|} s_{\mu/\nu}(x)s_{\mu/\nu}(y)
= \exp \sum_{n=1}^{\infty} \frac{u^n}{n(1 - (uv)^n)} p_n(x)p_n(y), \label{eqn:Schur2Sum} \\
&& \sum_{\mu, \nu} (-u)^{|\mu|} (-v)^{|\nu|} s_{\mu/\nu^t}(x)s_{\mu^t/\nu}(y)
= \exp \sum_{n=1}^{\infty} \frac{-u^n}{n(1 - (uv)^n)} p_n(x)p_n(y). \label{eqn:Schur2SumT}
\end{eqnarray}
For $|u_1u_2v_1v_2| < 1$, we have
\begin{equation}  \label{eqn:Schur4Sum}
\begin{split}
&  \sum_{\mu, \nu, \xi, \eta } u_1^{|\mu|} u_2^{|\nu|} v_1^{|\xi|} v_2^{|\eta|}
s_{\mu/\xi}(x) s_{\nu/\xi}(y) s_{\nu/\eta}(z)s_{\mu/\eta}(w) \\
= & \exp \sum_{n=1}^{\infty} \frac{1}{n(1- (u_1u_2v_1v_2)^n)}
(u_1^np_n(x)p_n(z) + u_2^n p_n(y)p_n(w) \\
& + (u_1u_2v_1)^n p_n(w)p_n(z) + (u_1u_2v_2)^n p_n(x)p_n(y)),
\end{split}
\end{equation}
and
\begin{equation}  \label{eqn:Schur4SumT}
\begin{split}
&  \sum_{\mu, \nu, \xi, \eta } (-u_1)^{|\mu|} (-u_2)^{|\nu|} v_1^{|\xi|} v_2^{|\eta|}
s_{\mu/\xi^t}(x) s_{\nu/\xi^t}(y) s_{\nu^t/\eta}(z)s_{\mu^t/\eta}(w) \\
= & \exp \sum_{n=1}^{\infty} \frac{1}{n(1- (u_1u_2v_1v_2)^n)}
(-u_1^np_n(x)p_n(z) - u_2^n p_n(y)p_n(w) \\
& + (u_1u_2v_1)^n p_n(w)p_n(z) + (u_1u_2v_2)^n p_n(x)p_n(y)),
\end{split}
\end{equation}
\end{lemma}

\begin{proof}
Write the left-hand side of (\ref{eqn:Schur2Sum}) as $S(x, y, u, v)$.
By (\ref{eqn:SchurSumT}) one has:
\begin{eqnarray*}
S(x, y; u, v)
& = & \sum_{\nu} (u v)^{|\nu|} \sum_{\mu} s_{\mu/\nu}(u x)s_{\mu/\nu}(y) \\
& = & \exp \sum_{n=1}^{\infty} \frac{u^n}{n} p_n(x)p_n(y) \cdot
\sum_{\nu, \mu} (u v)^{|\nu|} s_{\nu/\mu}(u x)s_{\nu/\mu}(y) \\
& = & \exp \sum_{n=1}^{\infty} \frac{u^n}{n} p_n(x)p_n(y) \cdot
\sum_{\nu, \mu} v^{|\nu|} u^{|\mu|} s_{\nu/\mu}(v y) s_{\nu/\mu}(u x) \\
& = &  \exp \sum_{n=1}^{\infty} \frac{u^n}{n} p_n(x)p_n(y) \cdot S(vy, ux; u, v).
\end{eqnarray*}
This procedure can be graphically represented as follows:
$$
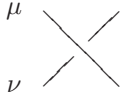
\xy
(0, 0); (10, -10), **@{-};
(0, -10); (4, -6), **@{-};
(6, -4); (10, 0), **@{-};
(-4, 0)*+{\mu}; (-4, -10)*+{\nu};
\endxy$$
One can repeat this procedure for infinitely many times to get:
\begin{eqnarray*}
&& \sum_{\mu, \nu} u^{|\mu|} v^{|\nu|} s_{\mu/\nu}(x)s_{\mu/\nu}(y) \\
& = & \exp \sum_{n=1}^{\infty} \frac{1+u^n+u^{2n} + (uv)^{2mn}}{n} p_n(x)p_n(y)
\cdot S((uv)^mx, (uv)^my; u) \\
& = & \exp \sum_{n=1}^{\infty} \frac{u^n (1 +(uv)^n+\cdots + (uv)^{2n} + \cdots)}{n} p_n(x)p_n(y) \\
& = & \exp \sum_{n=1}^{\infty} \frac{u^n}{n(1 - (uv)^n)} p_n(x)p_n(y).
\end{eqnarray*}
The proof of (\ref{eqn:Schur2SumT}) is similar.
We present two proof for (\ref{eqn:Schur4Sum}).
In the first proof we use (\ref{eqn:SchurSum}) and (\ref{eqn:SchurSumXY}) to reduce to
(\ref{eqn:Schur2Sum}):
\begin{eqnarray*}
&& \sum_{\mu, \nu, \xi, \eta} u_1^{|\mu|}u_2^{|\nu|}v_1^{|\xi|}v_2^{|\eta|}
s_{\mu/\xi}(x)s_{\nu/\xi}(y) s_{\mu/\eta}(z) s_{\nu/\eta}(w) \\
& = & \sum_{\mu, \xi, \eta} u_1^{|\mu|}u_2^{|\xi|}v_1^{|\xi|}v_2^{|\eta|}
s_{\mu/\xi}(x) s_{\mu/\eta}(z)
\sum_{\nu} s_{\nu/\xi}(u_2 y) s_{\nu/\eta}(w) \\
& = & \exp \sum_{n=1}^{\infty} \frac{u_2^n}{n} p_n(y)p_n(w) \cdot
\sum_{\mu, \nu, \xi, \eta} u_1^{|\mu|}u_2^{|\xi|}v_1^{|\xi|}v_2^{|\eta|}
s_{\mu/\xi}(x) s_{\mu/\eta}(z) s_{\eta/\nu}(u_2 y) s_{\xi/\nu}(w) \\
& = & \exp \sum_{n=1}^{\infty} \frac{u_2^n}{n} p_n(y)p_n(w) \cdot
\sum_{\mu, \nu} u_1^{|\mu|} (u_2v_1v_2)^{|\nu|} s_{\mu/\nu}(x, u_2v_1w) s_{\mu/\nu}(z, u_2v_2y) \\
& = & \exp \sum_{n=1}^{\infty} \frac{u_2^n}{n} p_n(y)p_n(w) \cdot
\exp \sum_{n=1}^{\infty} \frac{u_1^n}{n(1 - (u_1u_2v_1v_2)^n)} (p_n(x) \\
&& + (u_2v_1)^n p_n(w))(p_n(z) + (u_2v_2)^n p_n(y)) \\
& = & \exp \sum_{n=1}^{\infty} \frac{1}{n(1 - (u_1u_2v_1v_2)^n)} (u_1^n p_n(x) p_n(z) \\
&& + (u_1u_2v_2)^n p_n(x)p_n(y) + (u_1u_2v_1)^n p_n(w)p_n(z)
+ u_2^n p_n(y)p_n(w))
\end{eqnarray*}

The second proof of (\ref{eqn:Schur4Sum}) is similar to the proof of (\ref{eqn:Schur2Sum}). We write its left-hand side
as $S(x, y, z, w; u , v)$.
Then we have:
\begin{eqnarray*}
&& S(x, y, z, w; u_1, u_2, v_1, v_2) \\
& = & \sum_{\mu, \nu, \xi, \eta } u_1^{|\mu|}u_2^{|\nu|}v_1^{|\xi|}v_2^{|\eta|}
s_{\mu/\xi}(x) s_{\nu/\xi}(y) s_{\nu/\eta}(w) s_{\mu/\eta}(z) \\
& = &  \sum_{\xi, \eta } (u_1v_1)^{|\xi|} (u_2v_2)^{|\eta|}
\sum_{\mu} s_{\mu/\xi}(u_1 x) s_{\mu/\eta}(z)
\sum_{\nu} s_{\nu/\xi}(y) s_{\nu/\eta}(u_2w) \\
& = & \exp \sum_{n=1}^{\infty} \frac{1}{n} (u_1^np_n(x)p_n(z) + u_2^n p_n(y)p_n(w)) \\
&&  \cdot
\sum_{\mu, \nu, \xi, \eta} (u_1v_1)^{|\xi|} (u_2v_2)^{|\eta|}
s_{\eta/\mu}(u_1 x) s_{\xi/\mu}(z) s_{\eta/\nu}(y) s_{\xi/\nu}(u_2 w) \\
& = & \exp \sum_{n=1}^{\infty} \frac{1}{n} (u_1^np_n(x)p_n(z) + u_2^n p_n(y)p_n(w)) \\
&&  \cdot
\sum_{\mu, \nu, \xi, \eta} u_1^{|\xi|}v_1^{|\mu|} u_2^{|\eta|}v_2^{|\nu|}
s_{\eta/\mu}(u_1 x) s_{\xi/\mu}(v_1z) s_{\eta/\nu}(v_2y) s_{\xi/\nu}(u_2 w) \\
& = &  \exp \sum_{n=1}^{\infty} \frac{1}{n} (u_1^np_n(x)p_n(z) + u_2^n p_n(y)p_n(w))
\cdot S(v_1z, u_1x, u_2w, v_2y; u_1, u_2, v_1, v_2).
\end{eqnarray*}
Repeating this procedure for four times:
\begin{eqnarray*}
&& S(x, y, z, w; u_1, u_2, v_1, v_2) \\
& = &  \exp \sum_{n=1}^{\infty} \frac{1}{n}
(u_1^n (1 + (u_1u_2v_1v_2)^n) p_n(x)p_n(z)
+ u_2^n(1 + (u_1u_2v_1v_2)^n) p_n(y)p_n(w) \\
&& + (u_1u_2v_1)^n(1 + (u_1u_2v_1v_2)^n)p_n(w)p_n(z) \\
&& + (u_1u_2v_2)^n (1 +  (u_1u_2v_1v_2)^{n})  p_n(x)p_n(y)) \\
&& \cdot S(u_1u_2v_1v_2x, u_1u_2v_1v_2y, u_1u_2v_1v_2z, u_1u_2v_1v_2w; u_1, u_2, v_1, v_2).
\end{eqnarray*}
Repeating this procedure for infinitely many times:
\begin{eqnarray*}
&& S(x, y, z, w; u_1, u_2, v_1, v_2) \\
& = &  \exp \sum_{n=1}^{\infty} \frac{1}{n(1- (u_1u_2v_1v_2)^n)}
(u_1^np_n(x)p_n(z) + u_2^n p_n(y)p_n(w) \\
&& + (u_1u_2v_1)^n p_n(w)p_n(z) + (u_1u_2v_2)^n p_n(x)p_n(y)).
\end{eqnarray*}
\end{proof}

\subsection{Some results on specialization of Schur functions}

Recall any symmetric function $f$ can be written a polynomial
$f(e_1, e_2, \dots, e_n, \dots)$ in the elementary symmetric functions
$e_1, \dots, e_n, \dots$.
Let $E(u) = 1 + \sum_{n =1}^{\infty} e_n u^n$.
We write $f(e_1, \dots, e_n, \dots)$ as $f(E(u))$.
In general,
take $E(u)$ to be any formal power series with leading coefficient $1$,
$f(E(u))$ defines a specialization of $f$.

\begin{lemma} \label{lm:pn}
Let $a = (a_1, a_2, \dots)$, $b = (b_1, b_2, \dots)$,
and
$$E(u) = \frac{\prod_{i=1}^{\infty} (1 + a_i u)}{\prod_{i=1}^{\infty} (1 + b_i u)}.$$
Then we have
\begin{eqnarray} \label{eqn:pn}
&& p_n(E(u)) = p_n(a) - p_n(b)
\end{eqnarray}
\end{lemma}

\begin{proof}
Let
$$P(u) = \sum_{n=1}^{\infty} p_n(E(u)) u^{n-1}.$$
Then one has \cite{Mac}:
\begin{eqnarray*}
P(-u) & = & \frac{d}{du} \log E(u) = \sum_{i=1}^{\infty} (\frac{a_i}{1 + a_i u}
- \frac{b_i}{1 + b_i u} ) \\
& = & \sum_{i=1}^{\infty} \sum_{n=1}^{\infty} ( a_i^n  - b_i^n) (-u)^{n-1}
= \sum_{n=1}^{\infty} (p_n(a) - p_n(b)) (-u)^{n-1}.
\end{eqnarray*}
\end{proof}

\begin{lemma} \label{lm:Sum}
Let $a = (a_1, a_2, \dots)$, $b = (b_1, b_2, \dots)$,
$c = (c_1, c_2, \dots)$, $d = (d_1, \dots)$,
and
\begin{align*}
E(u) & = \frac{\prod_{i=1}^{\infty} (1 + a_i u)}{\prod_{i=1}^{\infty} (1 + b_i u)}, &
\tilde{E}(u) & = \frac{\prod_{i=1}^{\infty} (1 + c_i u)}{\prod_{i=1}^{\infty} (1 + d_i u)}.
\end{align*}
Then we have
\begin{eqnarray}
&& \sum_{\mu} s_\mu(E(u))s_\mu(\tilde{E}(u))
= \prod_{i, j =1}^{\infty} \frac{(1 - a_id_j)(1- b_ic_j)}{(1- a_ic_j)(1- b_i d_j)}.
\end{eqnarray}
\end{lemma}

\begin{proof}
This is a straightforward consequence of (\ref{eqn:SchurSum1}) and Lemma \ref{lm:pn}.
\end{proof}

\section{The Deformed Topological Vertex}
\label{sec:DTV}

In this section we will introduce the deformed topological vertex.
We will begin by studying the large $N$ Chern-Simons invariants of the Hopf link.
We will express them in terms of specializations of skew Schur functions
and study their symmetric properties.
It will be interesting to combine our approach with that of Awata and Konno \cite{Awa-Kon}.
We hope to report on this in a separate work.

\subsection{The quantum dimension}

Recall the large $N$ invariant of the unknot is given by the quantum dimension \cite{Mar-Vaf}:
\begin{eqnarray} \label{def:WMu}
&& W_{\mu}(q, e^{-t})
= \prod_{e \in \mu} \frac{e^{t/2}q^{c(e)/2} - e^{-t/2}q^{-c(e)/2}}{q^{h(e)/2}- q^{-h(e)/2}},
\end{eqnarray}
where $\mu = (\mu_1, \dots, \mu_l)$ is a partition.
For the zero partition we take $W_{(0)}(q, e^{-t}) = 1$.
Here $c(e)$ and $h(e)$ are the content and the hook length of the box $e$ in the Young diagram $\mu$ respectively.
From the definition (\ref{def:WMu}) one easily gets the following:

\begin{proposition}
The quantum dimension has the following symmetry properties:
\begin{eqnarray}
&& W_{\mu}(q^{-1}, e^{t})  = W_{\mu}(q, e^{-t}),  \label{eqn:WMut1} \\
&& W_{\mu^t}(q^{-1}, e^{-t}) = W_{\mu^t}(q, e^{t}) = (-1)^{|\mu|} W_{\mu}(q, e^{-t}). \label{eqn:WMut2}
\end{eqnarray}
\end{proposition}

Our next result express the quantum dimension as a specialization of the Schur functions.

\begin{proposition}
The quantum dimension can be identified with the following specialization
of the Schur function:
\begin{eqnarray} \label{eqn:WMu}
W_{\mu}(q, e^{-t}) = (qe^t)^{|\mu|/2} s_{\mu}(E^{(0)}(u; q, e^{-t})),
\end{eqnarray}
where
\begin{eqnarray} \label{eqn:E(u)}
&& E^{(0)}(u) = 1 + \sum_{n = 1}^{\infty} u^n e_n
= 1 + \sum_{n=1}^{\infty} u^n \prod_{i=1}^n \frac{1-e^{-t}q^{i-1}}{q^i - 1}.
\end{eqnarray}
\end{proposition}

\begin{proof}
Recall the following fact from the theory of symmetric functions.
For a specialization with
$$H(u) = \sum_{r =0}^{\infty} h_r u^r =  \prod_{i=0}^{\infty} \frac{1-bq^it}{1-aq^it},$$
one has (\cite{Mac}, p. 27, Example 5,  and p. 45, Example 3).
\begin{eqnarray}
h_r & = & \prod_{i=1}^r \frac{a-bq^{i-1}}{1-q^i}, \\
e_r & = & \prod_{i=1}^r \frac{aq^{i-1} - b}{1 - q^i}, \\
p_r & = & \frac{a^r - b^r}{1- q^r}, \\
s_{\nu} & = & q^{n(\nu)} \prod_{x \in \nu} \frac{a-bq^{c(x)}}{1 - q^{h(x)}}.
\end{eqnarray}
Hence one gets $E(u)$ in (\ref{eqn:E(u)}) by taking
$a = e^{-t}$, $b = 1$, and $q = e^{\sqrt{-1}\lambda}$.
It follows that
\begin{eqnarray*}
s_{\mu}(E(u; q, e^{-t})) & = & q^{n(\mu)}\prod_{e \in \mu} \frac{e^{-t} -  q^{c(e)}}{1 - q^{h(e)}} \\
& = & q^{n(\mu)+\frac{1}{2}\sum_{e\in \mu} [c(e) - h(e)]}e^{-|\mu|t/2}
\prod_{e \in \mu} \frac{e^{t/2}q^{c(e)/2} - e^{-t/2}q^{-c(e)/2}}{q^{h(e)/2} - q^{-h(e)/2}} \\
& = & (q^{-1} e^{-t})^{|\mu|/2}
\prod_{e \in \mu} \frac{e^{t/2}q^{c(e)/2} - e^{-t/2}q^{-c(e)/2}}{q^{h(e)/2}- q^{-h(e)/2}} \\
& = & (q^{-1}e^{-t})^{|\mu|/2} W_{\mu}(q, e^{-t}).
\end{eqnarray*}
Here we have used the following properties of the hook length and the content \cite{Mac}:
\begin{eqnarray}
&& \sum_{x \in \nu} h(x) = n(\nu) + n(\nu^t) + |\nu|, \\
&& \sum_{x \in \nu} c(x) = n(\nu^t) - n(\nu).
\end{eqnarray}
\end{proof}

\subsection{The invariant of the Hopf link}

Recall the large $N$ invariant of the Hopf link is \cite{Mor-Luk, Aga-Mar-Vaf}:
\begin{eqnarray} \label{def:WMuNu}
W_{\mu, \nu}(q, e^{-t}) = W_{\mu}(q, e^{-t}) (e^{t}q)^{|\nu|/2}
s_{\nu}(E^{\mu}(u; q, e^{-t})),
\end{eqnarray}
where
$$E^{\mu}(u; q, e^{-t})
= \prod_{i=1}^{l(\mu)} \frac{1+q^{\mu_i-i}u}{1+q^{-i}u} \cdot
(1 + \sum_{n=1}^{\infty} u^n \prod_{i=1}^n \frac{1-e^{-t}q^{i-1}}{q^i - 1}).$$

\begin{proposition}
We have
\begin{eqnarray} \label{eqn:EMu}
E^{\mu}(u) & = &
\frac{\prod_{i=1}^{\infty} (1 + q^{\mu_i-i}u)}{\prod_{i=1}^{\infty} (1 + e^{-t} q^{-i} u)}.
\end{eqnarray}
\end{proposition}

\begin{proof}
Recall the following famous identity due to Cauchy \cite{And}.
For $|x| < 1$, $|q| < 1$,
\begin{eqnarray*}
&& 1+ \sum_{n=1}^{\infty} x^n \prod_{i=1}^n \frac{(1 - a q^{i-1})}{( 1- q^i)}
= \prod_{i=0}^{\infty} \frac{(1- axq^i)}{(1- x q^{i})}.
\end{eqnarray*}
Now
\begin{eqnarray*}
&& 1 + \sum_{n=1}^{\infty} u^n \prod_{i=1}^n \frac{1-e^{-t}q^{i-1}}{q^i - 1}
= 1 + \sum_{n=1}^{\infty} (-e^{-t} q^{-1} u)^n \prod_{i=1}^n \frac{1-e^{t}q^{-(i-1)}}{1- q^{-i}} \\
& = & \prod_{i=1}^{\infty} \frac{1 + q^{-i}u}{1 + e^{-t} q^{-i}u}.
\end{eqnarray*}
Hence (\ref{eqn:EMu}) follows.
\end{proof}

Apply Lemma \ref{lm:pn} to $E^{\mu}(u)$:
\begin{eqnarray} \label{eqn:pnEMu}
p_n(E^{\mu}(u; q, e^{-t}))
& = & \sum_{i=1}^{\infty} (q^{n(\mu_i - i)}
- e^{-nt} q^{n(-i)}).
\end{eqnarray}
Now note
$$\sum_{i=1}^{\infty} q^{n(-i)} = \frac{q^{-n}}{1 - q^{-n}} = - \frac{1}{1-q^n}.$$
Hence we have
\begin{eqnarray} \label{eqn:pnEMu2}
p_n(E^{\mu}(u; q, e^{-t}))
& = & p_n(q^{\mu_i +\rho-\frac{1}{2}})+ p_n(e^{-t} q^{-\rho-\frac{1}{2}}).
\end{eqnarray}
I.e., formally $E^{\mu}(u; q, e^{-t})$ gives the specialization
$(q^{\mu+\rho-\frac{1}{2}}, e^{-t}q^{-\rho-\frac{1}{2}})$.
In (\ref{eqn:pnEMu2}) and the expression below,
we are not working in different domains $|q| < 1$ and $|q| > 1$ at the same time.
Instead,
we write the relevant skew Schur functions as polynomials of $p_n(q^{\rho})$ and $p_n(q^{-\rho})$,
whose values are given by (\ref{eqn:PnQrho}).
Hence our treatment is mathematically rigorous.
For later use note
\begin{eqnarray}
&& p_n(E^{\mu}(u; q, e^{-t})) = p_1(E^{\mu}(q^n, e^{-nt})),
\end{eqnarray}
and
\begin{eqnarray} \label{eqn:p1EMu}
&& p_1(E^{\mu}(u; q, e^{-t}))
= \sum_{i=1}^{\infty} (q^{\mu_i - i} - e^{-t} q^{-i})
= \sum_{i=1}^{l(\mu)} (q^{\mu_i - i} - q^{-i}) + \frac{1 - e^{-t}}{q- 1}.
\end{eqnarray}

\begin{lemma}
We have
\begin{eqnarray}
&& s_{\nu}(E^{\mu}(u; q, e^{-t}))
= q^{-|\nu|/2} \sum_{\eta} s_{\nu/\eta}(e^{-t}q^{-\rho}) s_{\eta}(q^{\mu+\rho})  \label{eqn:sNuEMu1} \\
& = & q^{-|\nu|/2} \sum_{\xi, \eta} s_{\nu/\eta}(e^{-t}q^{-\rho})  (-1)^{|\eta|} q^{\kappa_{\eta}/2}
\frac{s_{\mu/\xi}(q^{-\rho})s_{\eta/\xi}(q^{-\rho})}{s_{\mu}(q^{-\rho})}. \label{eqn:sNuEMu2}
\end{eqnarray}
\end{lemma}

\begin{proof}
We use (\ref{eqn:pnEMu2}) and (\ref{eqn:SchurSumXY}) to get:
\begin{eqnarray*}
s_{\nu}(E^{\mu}(u; q, e^{-t}))
& = & s_{\nu}(q^{\mu+\rho-\frac{1}{2}}, e^{-t}q^{-\rho-\frac{1}{2}})
= \sum_{\eta}  s_{\nu/\eta}(e^{-t}q^{-\rho-\frac{1}{2}})s_{\eta}(q^{\mu+\rho-\frac{1}{2}}) \\
& = & \sum_{\eta} q^{-(|\nu|-|\eta|)/2} s_{\nu/\eta}(e^{-t}q^{-\rho})
q^{-|\eta|/2} s_{\eta}(q^{\mu+\rho}) \\
& = &  q^{-|\nu|/2} \sum_{\eta} s_{\nu/\eta}(e^{-t}q^{-\rho}) s_{\eta}(q^{\mu+\rho}).
\end{eqnarray*}
This proves (\ref{eqn:sNuEMu1}).
To prove (\ref{eqn:sNuEMu2}) recall the following identity proved in \cite{ZhoConjecture}:
\begin{eqnarray}
&& s_{\nu}(q^{\mu + \rho}) = (-1)^{|\mu|} q^{\kappa_{\nu}/2}\sum_{\xi}
\frac{s_{\mu/\xi}(q^{-\rho})s_{\nu/\xi}(q^{-\rho})}{s_{\mu}(q^{-\rho})}.
\end{eqnarray}
\end{proof}

\begin{corollary}
We have the following symmetry:
\begin{eqnarray} \label{eqn:SymSNuEta}
s_{\nu^t/\eta^t}(E^{\mu^t}(u; q, e^{-t})) = (-1)^{|\nu|-|\eta|}
q^{-(|\nu|-|\eta|)} s_{\nu/\eta}(E^{\mu}(u; q^{-1}, e^{-t})).
\end{eqnarray}
\end{corollary}

\begin{proof}
Recall the following symmetry \cite{ZhoCounting}:
\begin{eqnarray}
s_{\lambda/\mu}(q^{\nu+\rho}) & = & (-1)^{|\lambda|-|\mu|}s_{\lambda^t/\mu^t}(q^{-\nu-\rho}).
\label{eqn:slambda/mu}
\end{eqnarray}
Hence (\ref{eqn:SymSNuEta}) follows from (\ref{eqn:sNuEMu2}).
\end{proof}

\subsection{Symmetries of $W_{\mu^1, \mu^2}$}

\begin{proposition}
We have
\begin{eqnarray}
&& W_{\mu^1, \mu^2}(q, e^{-t}) = W_{\mu^2, \mu^1}(q, e^{-t}), \label{eqn:Wmumu12} \\
&& W_{(\mu^1)^t, (\mu^2)^t}(q, e^{-t}) = (-1)^{|\mu^1|+|\mu^2|} W_{\mu^1, \mu^2}(q^{-1}, e^{-t}). \label{eqn:Wtt}
\end{eqnarray}
\end{proposition}

\begin{proof}
The identity (\ref{eqn:Wmumu12}) follows from the fact that $W_{\mu^1, \mu^2}(q, e^{-t})$ can be obtained from
the colored HOMFLY polynomials of the Hopf link.
By (\ref{def:WMuNu}) and (\ref{eqn:SymSN}),
\begin{eqnarray*}
W_{\mu^t, \nu^t}(q, e^{-t})
& = & W_{\mu^t}(q, e^{-t}) (e^{t}q)^{|\nu^t|/2}
s_{\nu^t}(E_{\mu^t}(u; q, e^{-t})) \\
& = & (-1)^{|\mu|} W_{\mu}(q^{-1}, e^{-t}) (e^{t}q)^{|\nu|/2}
(-1)^{|\nu|} q^{-|\nu|} s_{\nu}(E_{\mu}(u; q^{-1}, e^{-t})) \\
& = & (-1)^{|\mu|+|\nu|} W_{\mu, \nu}(q^{-1}, e^{-t}).
\end{eqnarray*}
This proves (\ref{eqn:Wtt}).
\end{proof}

\subsection{The Deformed topological vertex}

The topological vertex introduced in \cite{Aga-Kle-Mar-Vaf} is defined by
\begin{eqnarray} \label{eqn:TV}
\cW_{\mu^1, \mu^2, \mu^3}(q)
= \sum_{\rho^1, \rho^3}c_{\rho^1(\rho^3)^t}^{\mu^1(\mu^3)^t}q^{\kappa_{\mu^2}/2+\kappa_{\mu^3}/2}
\frac{\cW_{(\mu^2)^t\rho^1}(q)\cW_{\mu^2(\rho^3)^t}(q)}{\cW_{\mu^2}(q)},
\end{eqnarray}
where
$$ c_{\rho^1(\rho^3)^t}^{\mu^1(\mu^3)^t}
= \sum_{\eta} c_{\eta\rho^1}^{\mu^1}c_{\eta(\rho^3)^t}^{(\mu^3)^t}.
$$
It can be reformulated in terms of skew Schur functions \cite{ORV, ZhoCounting, Hol-Iqb-Vaf}:
\begin{eqnarray*}
\cW_{\mu^1, \mu^2, \mu^3}(q)
& = & q^{\kappa_{\mu^3}/2}
s_{\mu^2}(q^{\rho}) \sum_{\eta}
s_{\mu^1/\eta}(q^{(\mu^2)^t+\rho})
s_{(\mu^3)^t/\eta}(q^{\mu^2+\rho}).
\end{eqnarray*}
Now we define the deformed topological vertex to be:
\begin{equation} \label{def:W}
\begin{split}
& W_{\mu^1, \mu^2, \mu^3}(q, e^{-t}) \\
= & q^{\kappa_{\mu^3}/2} e^{-|\mu^2|t/2} W_{\mu^2}(q, e^{-t}) \\
& \cdot \sum_{\eta}
q^{(|\mu^1| - |\eta|)/2} s_{\mu^1/\eta}(E^{(\mu^2)^t}(u; q, e^{-t}))
q^{(|\mu^3| - |\eta|)/2} s_{(\mu^3)^t/\eta}(E^{\mu^2}(u; q, e^{-t})).
\end{split}
\end{equation}
It is easy to see that the leading term of $W_{\mu^1, \mu^2, \mu^3}(q, e^{-t})$
is $\cW_{\mu^1, \mu^2, \mu^3}(q)$.
When $\mu^1 = (0)$,
\begin{eqnarray*}
W_{(0), \mu^2, \mu^3}(q, e^{-t})
& = & q^{\kappa_{\mu^3}/2} e^{-|\mu^2|t/2}W_{\mu^2}(q, e^{-t})
q^{|\mu^3|/2} s_{(\mu^3)^t}(E^{\mu^2}(u; q, e^{-t})) \\
& = & e^{-(|\mu^2|+|\mu^3|)t/2} q^{\kappa_{\mu^3}/2} W_{\mu^2, (\mu^3)^t}(q, e^{-t}).
\end{eqnarray*}
When $\mu^3 = (0)$,
\begin{eqnarray*}
W_{\mu^1, \mu^2, (0)}(q, e^{-t})
& = & e^{-|\mu^2|t/2} W_{\mu^2}(q, e^{-t}) q^{|\mu^1|/2} s_{\mu^1}(E^{(\mu^2)^t}(u; q, e^{-t})).
\end{eqnarray*}
When $\mu^2 = (0)$,
\begin{equation} \label{eqn:WMu2=0}
\begin{split}
& W_{\mu^1, (0), \mu^3}(q, e^{-t}) \\
= & q^{\kappa_{\mu^3}/2}
\sum_{\eta}
q^{(|\mu^1| - |\eta|)/2} s_{\mu^1/\eta}(E(u; q, e^{-t}))
q^{(|\mu^3| - |\eta|)/2} s_{(\mu^3)^t/\eta}(E(u; q, e^{-t})).
\end{split}
\end{equation}
One then sees that $W_{\mu^1, \mu^2, \mu^3}$ does not have the $\bZ_3$ cyclic symmetry.

\begin{proposition}
We have
\begin{eqnarray}
&& W_{(\mu^1)^t, (\mu^2)^t, (\mu^3)^t}(q, e^{-t})
= (-1)^{|\mu^1|+|\mu^2|+|\mu^3|} W_{\mu^3, \mu^2, \mu^1}(q^{-1}, e^{-t}).
\end{eqnarray}
\end{proposition}

\begin{proof}
This is a straightforward consequence of (\ref{eqn:SymSNuEta}).
\end{proof}

\section{From Summations to Products}
\label{sec:Product}

We use the results in the preceding two sections to prove some product
expressions for certain sums of the invariants of the Hopf link and the deformed topological vertex.
They generalize the results in \cite{Iqb-Kas1, Iqb-Kas2, Egu-Kan1, Egu-Kan2, ZhoCounting}
They will play a crucial role in the next section.

\subsection{Results on $F_{\mu^1\mu^2}$}
\label{sec:fmumu}

For two partitions $\mu^1, \mu^2$, define:
\begin{eqnarray*}
F_{\mu^1 \mu^2}(q, e^{-t})
& = & \sum_{i \geq 1} (q^{\mu^1_i-i+\frac{1}{2}} - e^{-t} q^{-i+\frac{1}{2}})
\sum_{j \geq 1} (q^{\mu^2_j-j+\frac{1}{2}} - e^{-t} q^{-j+\frac{1}{2}}) \\
&& - \sum_{i \geq 1} (1- e^{-t}) q^{-i+\frac{1}{2}} \sum_{j \geq 1}(1- e^{-t}) q^{-j+\frac{1}{2}}.
\end{eqnarray*}
We have
\begin{eqnarray} \label{eqn:Ff}
&& F_{\mu^1\mu^2}(q, e^{-t}) = f_{\mu^1, \mu^2}(q) - e^{-t} f_{\mu^1}(q) - e^{-t}f_{\mu^2}(q).
\end{eqnarray}
Indeed, by (\ref{eqn:p1EMu}) we have
\begin{eqnarray*}
&& F_{\mu^1 \mu^2}(q, e^{-t}) \\
& = & q p_1(E^{\mu^1}(u; q, e^{-t})) p_1(E^{\mu^2}(u; q, e^{-t})
- q p_1(E^{(0)}(u; q, e^{-t})) p_1(E^{(0)}(u; q, e^{-t}) \\
& = & q \left(\frac{1 - e^{-t}}{q-1} + \sum_{i = 1}^{l(\mu^1)} (q^{\mu^1_i-i} -q^{-i})\right)
\left(\frac{1 - e^{-t}}{q-1} + \sum_{j = 1}^{l(\mu^2)} (q^{\mu^2_j-j} -q^{-j})\right)
- \frac{(1- e^{-t})^2q}{(1-q)^2} \\
& = & (1-e^{-t}) \frac{q}{q-1}\sum_{i = 1}^{l(\mu^1)} (q^{\mu^1_i-i} -q^{-i})
+ (1 - e^{-t}) \frac{q}{q-1} \sum_{j = 1}^{l(\mu^2)} (q^{\mu^2_j-j} -q^{-j}) \\
&& + q \sum_{i = 1}^{l(\mu^1)} (q^{\mu^1_i-i} -q^{-i}) \cdot \sum_{j = 1}^{l(\mu^2)} (q^{\mu^2_j-j} -q^{-j}) \\
& = & (1 - e^{-t}) f_{\mu^1}(q) + (1 - e^{-t}) f_{\mu^2} (q) + (q - 2 + q^{-1}) f_{\mu^1}(q) f_{\mu^2}(q) \\
& = & f_{\mu^1\mu^2}(q) - e^{-t} f_{\mu^1}(q) - e^{-t}f_{\mu^2}(q).
\end{eqnarray*}
The expression
\begin{eqnarray} \label{def:fMuNu}
&& f_{\mu^1\mu^2} = \sum_{i , j=1}^{\infty} (q^{\mu^1_i + \mu^2_j - i - j + 1} - q^{-i-j+1})
\end{eqnarray}
has been studied in \cite{ZhoCounting}:
\begin{eqnarray}
f_{\mu^1(\mu^2)^t}(q)
& = & \sum_{(i,j)\in \mu^1} q^{\mu^1_i+(\mu^2)^t_j-i-j+1}
+ \sum_{(i,j) \in \mu^2} q^{-(\mu^1)^t_j-\mu^2_i+i+j-1}  \label{eqn:fMuNu1} \\
& = & - \sum_{i, j \geq 1} (q^{\mu^1_i-\mu^2_j+j-i} - q^{j-i}).  \label{eqn:fMuNu2}
\end{eqnarray}
In particular,
\begin{eqnarray}
&& f_{\mu}(q) = \sum_{(i,j) \in \mu} q^{j-i}. \label{eqn:fMu}
\end{eqnarray}

\subsection{Some product expressions}

By (\ref{eqn:Ff}) we have
\begin{eqnarray*}
&& \exp \left(\sum_{n \geq 1} \frac{u^n}{n} F_{\mu^1(\mu^2)^t}(q^n, e^{-nt})\right) \\
& = & \exp \left( \sum_{n \geq 1} \frac{u^n}{n}
(f_{\mu^1(\mu^2)^t}(q^n) - e^{-nt} f_{\mu^1}(q^n) - e^{-nt}f_{(\mu^2)^t}(q^n))\right).
\end{eqnarray*}
From this one can obtain various product expressions.
First of all,
by (\ref{def:fMuNu}),
\begin{equation} \label{eqn:FProd1}
\begin{split}
& \exp \left(\sum_{n \geq 1} \frac{u^n}{n} F_{\mu^1(\mu^2)^t}(q^n, e^{-nt})\right) \\
= & \prod_{i, j =1}^{\infty}
\frac{(1- u q^{- i - j + 1})(1- u e^{-t} q^{\mu^1_i - i - j + 1})(1- u e^{-t} q^{(\mu^2)^t_j - i - j + 1})}
{(1- u q^{\mu^1_i + (\mu^2)^t_j - i - j + 1})(1- u e^{-t} q^{- i - j + 1})(1- u e^{-t}q^{- i - j + 1})}.
\end{split}
\end{equation}
By (\ref{eqn:fMuNu2}),
\begin{equation}  \label{eqn:FProd2}
\begin{split}
& \exp \left(\sum_{n \geq 1} \frac{u^n}{n} F_{\mu^1(\mu^2)^t}(q^n, e^{-nt})\right) \\
= & \prod_{i, j =1}^{\infty}
\frac{(1- u q^{\mu^1_i - \mu^2_j - i + j})(1- u e^{-t} q^{j - i})(1- u e^{-t}q^{j - i})}
{(1- u q^{j - i})(1- u e^{-t} q^{\mu^1_i - i + j})(1- u e^{-t} q^{\mu^2_j - i + j})}.
\end{split} \end{equation}
By (\ref{eqn:fMuNu1}) and (\ref{eqn:fMu}),
\begin{equation}  \label{eqn:FProd3}
\begin{split}
& \exp \left(\sum_{n \geq 1} \frac{u^n}{n} F_{\mu^1(\mu^2)^t}(q^n, e^{-nt})\right) \\
= & \prod_{(i, j) \in \mu^1}
\frac{1 - ue^{-t} q^{j - i}}{1 - uq^{\mu^1_i + (\mu^2)^t_j - i - j + 1}}
\cdot \prod_{(i, j) \in \mu^2}
\frac{1 - ue^{-t} q^{-(j - i)}}{1 - uq^{-[(\mu^1)^t_j + \mu^2_i - i - j + 1]}}.
\end{split}
\end{equation}

\subsection{Product expressions for $\cK_{\mu^1(\mu^2)^t}(Q)$}

\label{sec:Kmumu}

Define
\begin{eqnarray*}
&& \cK_{\mu^1\mu^2}(Q; q, e^{-t})
= \sum_{\nu} Q^{|\nu|} e^{- (|\mu^1| + |\mu^2| + 2 |\nu|)t/2}
W_{\mu^1\nu}(q, e^{-t})W_{\nu\mu^2}(q, e^{-t}).
\end{eqnarray*}

\begin{proposition} \label{prop:K}
The following identities holds:
\begin{equation}
\begin{split}
& \frac{\cK_{\mu^1\mu^2}(Q; q, e^{-t})}{\cK_{(0)(0)}(Q; q, e^{-t})} \\
= & e^{-(|\mu^1|+|\mu^2|)t/2}W_{\mu^1}(q, e^{-t}) W_{\mu^2}(q, e^{-t})
\exp \left( \sum_{n=1}^{\infty} \frac{Q^n}{n}F_{\mu^1\mu^2}(q^n, e^{-nt})\right).
\end{split}
\end{equation}
\end{proposition}

\begin{proof}
We have
\begin{eqnarray*}
&& \cK_{\mu^1\mu^2}(Q) \\
& = & \sum_{\nu} Q^{|\nu|} e^{-(|\mu^1|+|\mu^2|)t/2}W_{\mu^1} q^{|\nu|/2} s_{\nu}(E^{\mu^1}(u;q, e^{-t}))
W_{\mu^2} q^{|\nu|/2} s_{\nu}(E^{\mu^2}(u; q, e^{-t})) \\
& = & e^{-(|\mu^1|+|\mu^2|)t/2}W_{\mu^1}  W_{\mu^2}
\exp \sum_{n =1}^{\infty} \frac{(qQ)^n}{n}
p_n(E^{\mu^1}(u; q, e^{-t})) p_n(E^{\mu^2}(u; q, e^{-t})).
\end{eqnarray*}
In particular,
when $\mu^1=\mu^2=(0)$,
\begin{eqnarray} \label{eqn:K00}
\cK_{(0)(0)}(Q)
& = & \exp \sum_{n =1}^{\infty} \frac{(qQ)^n}{n} p_n(E^{(0)}(u; q, e^{-t}))^2.
\end{eqnarray}
Hence
\begin{eqnarray*}
&& \frac{\cK_{\mu^1\mu^2}(Q)}{\cK_{(0)(0)}(Q)} \\
& = & e^{-(|\mu^1|+|\mu^2|)t/2} W_{\mu^1}  W_{\mu^2} \\
&& \exp \sum_{n =1}^{\infty} \frac{(qQ)^n}{n}
(p_n(E^{\mu^1}(u; q, e^{-t})) p_n(E^{\mu^2}(u; q, e^{-t})) - p_n(E^{(0)}(u; q, e^{-t}))^2)\\
& = & e^{-(|\mu^1|+|\mu^2|)t/2}W_{\mu^1}  W_{\mu^2}
\exp \left(\sum_{n \geq 1} \frac{(e^tQ)^n}{n} F_{\mu^1\mu^2}(q^n, e^{-nt})\right).
\end{eqnarray*}
\end{proof}

We will find product expressions for $\cK_{\mu^1(\nu^2)^t}$ in two different ways.
By (\ref{eqn:K00}), (\ref{eqn:EMu}) and Lemma \ref{lm:Sum} we have
\begin{eqnarray}
\cK_{(0)(0)}
& = &  \prod_{i, j=1}^{\infty}
\frac{(1-q^{-(i+j-1)}Q)^2}{(1- q^{-(i+j-1)}e^tQ)(1-q^{-(i+j-1)}e^{-t}Q)} \\
& = & \prod_{n=1}^{\infty} \frac{(1- q^{-n}Q)^2}{(1- q^{-n}e^tQ)(1-q^{-n}e^{-t}Q)}.
\end{eqnarray}
Recall
\begin{eqnarray}
&& W_{\mu}(q, e^{-t}) = \prod_{(i, j) \in \mu}
\frac{e^{t/2} q^{(j-i)/2} - e^{-t/2} q^{-(j-i)/2}}
{q^{(\mu_i + \mu^t_j - i - j + 1)/2} - q^{- (\mu_i + \mu^t_j - i - j +1)/2}}.
\end{eqnarray}
It is possible to rewrite this as an infinite product.
For the denominator, recall \cite[Lemma 2.1]{ZhoCounting},
\begin{eqnarray}
&& \sum_{1 \leq i < j < \infty} (t^{\mu_i-\mu_j + j - i} - t^{j-i})
= - \sum_{(i, j) \in \mu} t^{\mu_i + \mu^t_j - i - j + 1};
\end{eqnarray}
for the numerator, recall
\begin{eqnarray}
&& f_{\mu}(q) = \sum_{(i, j) \in \mu} q^{j-i}
= -\sum_{i, j = 1}^{\infty} (q^{\mu_i + j - i} - q^{j-i}).
\end{eqnarray}
Hence
\begin{eqnarray}
W_{\mu}(q, e^{-t})
& = & \prod_{1 \leq i <  j < \infty}
\frac{q^{(\mu_i - \mu_j - i + j)/2} - q^{- (\mu_i - \mu_j - i + j)/2}}{q^{(j- i)/2} - q^{- (j - i)/2}} \\
&& \cdot \prod_{i, j = 1}^{\infty}
\frac{e^{t/2}q^{(\mu_i + j - i)/2} - e^{-t/2} q^{- (\mu_i + j - i)/2}}
{e^{t/2} q^{(j-i)/2} - e^{-t/2} q^{-(j-i)/2}}. \nonumber
\end{eqnarray}

By (\ref{eqn:FProd2}) and (\ref{eqn:FProd3}) one then gets two product expressions as follows.

\begin{theorem} \label{thm:K}
We have
\begin{equation}
\begin{split}
& \frac{\cK_{\mu^1(\mu^2)^t}(Q)}{\cK_{(0)(0)}(Q)} \\
= & e^{- (|\mu^1| + |\mu^2|)t/2} \prod_{(i, j) \in \mu^1}
\frac{e^{t/2} q^{(j-i)/2} - e^{-t/2} q^{-(j-i)/2}}
{q^{(\mu^1_i + (\mu^1)^t_j - i - j + 1)/2} - q^{- (\mu^1_i + (\mu^1)^t_j - i - j +1)/2}} \\
& \cdot \prod_{(i, j) \in (\mu^2)^t}
\frac{e^{t/2} q^{(j-i)/2} - e^{-t/2} q^{-(j-i)/2}}
{q^{((\mu^2)^t_i + \mu^2_j - i - j + 1)/2} - q^{- ((\mu^2)^t_i + \mu^2_j - i - j +1)/2}} \\
& \cdot \prod_{(i, j) \in \mu^1}
\frac{1 - e^{-t}Qq^{j - i}}{1 - Qq^{\mu^1_i + (\mu^2)^t_j - i - j + 1}}
\cdot \prod_{(i, j) \in \mu^2}
\frac{1 - e^{-t} Q q^{-(j - i)}}{1 - Qq^{-[(\mu^1)^t_j + \mu^2_i - i - j + 1]}}.
\end{split}
\end{equation}
Furthermore,
\begin{equation}
\begin{split}
&
 \frac{\cK_{\mu^1(\mu^2)^t}(Q)}{\cK_{(0)(0)}(Q)} \\
= & e^{- (|\mu^1| + |\mu^2| + 2 |\nu|)t/2} \prod_{1 \leq i <  j < \infty}
\frac{q^{(\mu^1_i - \mu^1_j - i + j)/2} - q^{- (\mu^1_i - \mu^1_j - i + j)/2}}{q^{(j- i)/2} - q^{- (j - i)/2}} \\
& \cdot \prod_{i, j = 1}^{\infty}
\frac{e^{t/2}q^{(\mu^1_i + j - i)/2} - e^{-t/2} q^{- (\mu^1_i + j - i)/2}}
{e^{t/2} q^{(j-i)/2} - e^{-t/2} q^{-(j-i)/2}} \\
& \cdot \prod_{1 \leq i <  j < \infty}
\frac{q^{(\mu^2_i - \mu^2_j - i + j)/2} - q^{- (\mu^2_i - \mu^2_j - i + j)/2}}{q^{(j- i)/2} - q^{- (j - i)/2}} \\
& \cdot \prod_{i, j = 1}^{\infty}
\frac{e^{t/2}q^{(\mu^2_i + j - i)/2} - e^{-t/2} q^{- (\mu^2_i + j - i)/2}}
{e^{t/2} q^{(j-i)/2} - e^{-t/2} q^{-(j-i)/2}} \\
& \cdot \prod_{i, j =1}^{\infty}
\frac{(1- u q^{\mu^1_i - \mu^2_j - i + j})(1- u e^{-t} q^{j - i})(1- u e^{-t}q^{j - i})}{(1- u q^{j - i})(1- u e^{-t} q^{\mu^1_i - i + j})(1- u e^{-t} q^{\mu^2_j - i + j})}.
\end{split}
\end{equation}
\end{theorem}

\subsection{Product expressions for $\cL_{\mu^1(\mu^2)^t}(Q)$}

\label{sec:Lmumu}

Define
\begin{eqnarray*}
&& \cL_{\mu^1\mu^2}(Q; q, e^{-t})
= \sum_{\nu} e^{-(|\mu^1|+|\mu^2|+2|\nu|)t/2}
(-Q)^{|\nu|} W_{\mu^1\nu}(q, e^{-t})W_{\nu^t\mu^2}(q, e^{-t}).
\end{eqnarray*}

\begin{proposition} \label{prop:L}
The following identity holds:
\begin{equation}
\begin{split}
& \frac{\cL_{\mu^1\mu^2}(Q; q, e^{-t})}{\cL_{(0)(0)}(Q; q, e^{-t})} \\
= & e^{-(|\mu^1|+|\mu^2|)t/2} W_{\mu^1}(q, e^{-t}) W_{\mu^2}(q, e^{-t})
\exp \left( -\sum_{n=1}^{\infty} \frac{Q^n}{n}F_{\mu^1\mu^2}(q^n, e^{-nt})\right).
\end{split}
\end{equation}
\end{proposition}

\begin{proof}
We have
\begin{eqnarray*}
&& \cL_{\mu^1\mu^2}(Q) \\
& = & \sum_{\nu} (-Q)^{|\nu|} e^{-|\mu^1|t/2} W_{\mu^1} q^{|\nu|/2} s_{\nu}(E^{\mu^1}(u;q, e^{-t})) \\
&& \cdot e^{-|\mu^2|t/2}W_{\mu^2} q^{|\nu|/2} s_{\nu^t}(E^{\mu^2}(u; q, e^{-t})) \\
& = & e^{-(|\mu^1|+|\mu^2|)t/2} W_{\mu^1}  W_{\mu^2}
\exp \left(- \sum_{n =1}^{\infty} \frac{(qQ)^n}{n}
p_n(E^{\mu^1}(u; q, e^{-t})) p_n(E^{\mu^2}(u; q, e^{-t}))\right).
\end{eqnarray*}
In particular,
when $\mu^1=\mu^2=(0)$,
\begin{eqnarray} \label{eqn:L00}
\cL_{(0)(0)}(Q)
& = & \exp \left(-\sum_{n =1}^{\infty} \frac{(qQ)^n}{n} p_n(E^{(0)}(u; q, e^{-t}))^2)\right).
\end{eqnarray}
Hence
\begin{eqnarray*}
&& \frac{\cL_{\mu^1\mu^2}(Q)}{\cL_{(0)(0)}(Q)} \\
& = & e^{-(|\mu^1|+|\mu^2|)t/2} W_{\mu^1}  W_{\mu^2}
\exp \left(- \sum_{n =1}^{\infty} \frac{(qQ)^n}{n}
(p_n(E^{\mu^1}(u; q, e^{-t})) p_n(E^{\mu^2}(u; q, e^{-t})) \right. \\
&& \left. - p_n(E^{(0)}(u; q, e^{-t}))^2)\right)\\
& = & e^{-(|\mu^1|+|\mu^2|)t/2} W_{\mu^1}  W_{\mu^2}
\exp \left(-\sum_{n \geq 1} \frac{Q^n}{n} F_{\mu^1\mu^2}(q^n, e^{-nt})\right).
\end{eqnarray*}
\end{proof}

By (\ref{eqn:FProd2}) and (\ref{eqn:FProd3}) one then gets two product expressions.

\begin{theorem}
\begin{eqnarray*}
&& \frac{\cL_{\mu^1(\mu^2)^t}(Q)}{\cL_{(0)(0)}(Q)} \\
& = & \prod_{(i, j) \in \mu^1}
\frac{e^{t/2} q^{(j-i)/2} - e^{-t/2} q^{-(j-i)/2}}
{q^{(\mu^1_i + (\mu^1)^t_j - i - j + 1)/2} - q^{- (\mu^1_i + (\mu^1)^t_j - i - j +1)/2}} \\
&& \cdot \prod_{(i, j) \in (\mu^2)^t}
\frac{e^{t/2} q^{(j-i)/2} - e^{-t/2} q^{-(j-i)/2}}
{q^{((\mu^2)^t_i + \mu^2_j - i - j + 1)/2} - q^{- ((\mu^2)^t_i + \mu^2_j - i - j +1)/2}} \\
&& \cdot \prod_{(i, j) \in \mu^1}
\frac{1 - Qq^{\mu^1_i + (\mu^2)^t_j - i - j + 1}}{1 - Qe^{-t}q^{j - i}}
\cdot \prod_{(i, j) \in \mu^2}
\frac{1 - Qq^{-((\mu^1)^t_j + \mu^2_i - i - j + 1)}}{1 - Q e^{-t} q^{-(j - i)}}.
\end{eqnarray*}
In particular, when $\mu^1 = \mu^2 = \mu$,
\begin{equation} \label{eqn:LMuMu}
\begin{split}
& \frac{\cL_{\mu\mu^t}(Q)}{\cL_{(0)(0)}(Q)} \\
= & (-1)^{|\mu|} \prod_{e \in \mu}
\frac{(1 - e^{-t} q^{c(e)})(1 - e^{-t} q^{-c(e)})}{(1- Qe^{-t}q^{c(e)})(1-Qe^{-t}q^{-c(e)})}
\cdot \frac{(1 - Qq^{h(e)})(1 - Qq^{-h(e)})}{(1 - q^{h(e)})(1 - q^{-h(e)})}.
\end{split}
\end{equation}
Furthermore,
\begin{equation}
\begin{split}
& \frac{\cL_{\mu^1(\mu^2)^t}(Q)}{\cL_{(0)(0)}(Q)} \\
= & \prod_{1 \leq i <  j < \infty}
\frac{q^{(\mu^1_i - \mu^1_j - i + j)/2} - q^{- (\mu^1_i - \mu^1_j - i + j)/2}}{q^{(j- i)/2} - q^{- (j - i)/2}} \\
& \cdot \prod_{i, j = 1}^{\infty}
\frac{e^{t/2}q^{(\mu^1_i + j - i)/2} - e^{-t/2} q^{- (\mu^1_i + j - i)/2}}
{e^{t/2} q^{(j-i)/2} - e^{-t/2} q^{-(j-i)/2}} \\
& \cdot \prod_{1 \leq i <  j < \infty}
\frac{q^{(\mu^2_i - \mu^2_j - i + j)/2} - q^{- (\mu^2_i - \mu^2_j - i + j)/2}}{q^{(j- i)/2} - q^{- (j - i)/2}} \\
& \cdot \prod_{i, j = 1}^{\infty}
\frac{e^{t/2}q^{(\mu^2_i + j - i)/2} - e^{-t/2} q^{- (\mu^2_i + j - i)/2}}
{e^{t/2} q^{(j-i)/2} - e^{-t/2} q^{-(j-i)/2}} \\
& \cdot \prod_{i, j =1}^{\infty}
\frac{(1- u q^{\mu^1_i - \mu^2_j - i + j})(1- u e^{-t} q^{j - i})(1- u e^{-t}q^{j - i})}
{(1- u q^{j - i})(1- u e^{-t} q^{\mu^1_i - i + j})(1- u e^{-t} q^{\mu^2_j - i + j})}.
\end{split}
\end{equation}
\end{theorem}

\subsection{Generalizations}

Define:
\begin{eqnarray*}
&& \tilde{\cK}_{\mu^1\cdots \mu^N}(Q_{1}, \dots, Q_{N-1}; q, e^{-t})
= \sum_{\nu^1, \dots, \nu^{N-1}} \prod_{k=1}^{N}
q^{\kappa_{\nu^k}/2} W_{\nu^{k-1}\mu^k(\nu^{k})^t}(q, e^{-t}) Q_{k}^{|\nu^k|},
\end{eqnarray*}
where $\nu^0=\nu^N = (0)$, $Q_N = 1$.

\begin{proposition} \label{prop:KGen}
We have the following identity:
\begin{eqnarray*}
&& \frac{\tilde{\cK}_{\mu^1\cdots \mu^N}(Q_{1}, \dots, Q_{N-1})}
{\prod_{1 \leq k < l \leq N} \cK_{(0)(0)}(Q_k\cdots Q_{l-1})} \\
& = & \prod_{k=1}^{N} W_{\mu^k}(q, e^{-t}) \cdot
\exp \left( \sum_{1 \leq k < l \leq N}
\sum_{n \geq 1} \frac{(Q_k\cdots Q_{l-1})^n}{n} F_{\mu^i(\mu^j)^t}(q^n, e^{-nt}) \right).
\end{eqnarray*}
\end{proposition}

\begin{proof}
\begin{eqnarray*}
&& \tilde{\cK}_{\mu^1\cdots \mu^N}(Q_{1}, \dots, Q_{N-1}; q, e^{-t}) \\
& = & \sum_{\nu^1, \dots, \nu^{N-1}} \prod_{k=1}^{N}
q^{\kappa_{\nu^k}/2} W_{\nu^{i-1}\mu^k(\nu^{k})^t}(q, e^{-t}) Q_{k}^{|\nu^k|}
\hspace{.2in} (\nu^0=\nu^N = (0)) \\
& = & \sum_{\nu^1, \dots, \nu^{N-1}}
\prod_{k=1}^{N} e^{-|\mu^k|t/2} W_{\mu^k}(q, e^{-t}) \\
&& \cdot \sum_{\eta^{k-1}} q^{(|\nu^{k-1}|- |\eta^{k-1}|)/2} s_{\nu^{k-1}/\eta^{k-1}}(E^{(\mu^k)^t}q, e^{-t})) \\
&& \cdot q^{(|\nu^{k}|- |\eta^{k-1}|)/2}
s_{\nu^k/\eta^{k-1}}(E^{\mu^k}(q, e^{-t})) Q_k^{|\nu^k|} \\
&& \hspace{2in} (\nu^0=\nu^N = \eta^0=\eta^{N-1} = (0)) \\
& = & \prod_{k=1}^{N} e^{-|\mu^k|t/2}W_{\mu^k}(q, e^{-t})
\sum_{\nu^1, \dots, \nu^{N-1}}\sum_{\eta^1, \dots, \eta^{N-2}} \prod_{k=1}^{N-1}
q^{(|\nu^{k}|- |\eta^{k-1}|)/2} \\
&& \cdot s_{\nu^k/\eta^{k-1}}(E^{\mu^k}(q, e^{-t})) q^{(|\nu^{k}|- |\eta^{k}|)/2} Q_k^{|\nu^k|}
s_{\nu^k/\eta^k}(E^{(\mu^{k+1})^t}(q, e^{-t})) \\
&& \hspace{2in} (\eta^0=\eta^{N-1} = (0)) \\
& = &  \prod_{k=1}^{N} e^{-|\mu^k|t/2}W_{\mu^k}(q, e^{-t})) \cdot
\prod_{1 \leq k < l \leq N} \exp \sum_{n=1}^{\infty}
\frac{(qQ_kQ_{k+1} \cdots Q_{l})^n}{n} \\
&& \cdot p_n(E^{\mu^k}(q, e^{-t}))p_n(E^{(\mu^l)^t}(q, e^{-t})).
\end{eqnarray*}
In the last equality we have used (\ref{eqn:SchurSum3}).
Hence
\begin{eqnarray*}
&& \frac{\tilde{\cK}_{\mu^1\cdots \mu^N}(Q_{1}, \dots, Q_{N-1}; q, e^{-t})}
{\prod_{1 \leq k < l \leq N} \cK_{(0)(0)}(Q_k\cdots Q_{l-1}; q, e^{-t})} \\
& = & \prod_{k=1}^{N} (e^{-|\mu^k|t/2} W_{\mu^k}(q, e^{-t})) \cdot
\exp \sum_{1 \leq k < l \leq N}
\sum_{n \geq 1} \frac{(Q_k\cdots Q_{l-1})^n}{n} F_{\mu^k(\mu^l
)^t}(q^n, e^{-nt}).
\end{eqnarray*}
\end{proof}

Hence one can easily get a Nekrasov type infinite product expression and
the following product expression:
\begin{equation} \label{eqn:TildeK}
\begin{split}
& \frac{\tilde{\cK}_{\mu^1\cdots \mu^N}(Q_{1}, \dots, Q_{N-1}; q, e^{-t})}
{\prod_{1 \leq k < l \leq N} \cK_{(0)(0)}(Q_k\cdots Q_{l-1}; q, e^{-t})} \\
= & \prod_{\alpha=1}^N (e^{-|\mu^{\alpha}|t/2} \prod_{(i, j) \in \mu^{\alpha}}
\frac{e^{t/2} q^{(j-i)/2} - e^{-t/2} q^{-(j-i)/2}}
{q^{h(i,j)/2} - q^{-h(i,j)/2}}) \\
& \cdot \prod_{1 \leq \alpha < \gamma \leq N}
( \prod_{(i,j)\in \mu^{\alpha}} \frac{1 - Q_{\alpha} \cdots Q_{\gamma-1}e^{-t}q^{j-i}}
{1- Q_{\alpha} \cdots Q_{\gamma-1} q^{\mu^{\alpha}_i + (\mu^{\gamma})^t_j - i + j + 1}} \\
& \cdot  \prod_{(i,j)\in \mu^{\gamma}}  \frac{1 -  Q_{\alpha} \cdots Q_{\gamma-1}e^{-t}q^{-(j-i)}}
{1- Q_{\alpha} \cdots Q_{\gamma-1} q^{\mu^{\alpha}_j + (\mu^{\gamma})^t_i - i + j + 1}})
\end{split}
\end{equation}

\section{Deformed Topological String Partition Functions of Some Local Calabi-Yau Geometries}
\label{sec:CY}

In this section we introduce some deformed topological string
partition functions for some local Calabi-Yau geometries.
We identify them with the deformed instanton counting partition functions considered in
\S \ref{sec:Instanton}.
Our results generalize those related to Nekrasov's partition function \cite{Nek, Iqb-Kas1, Egu-Kan1, ZhoCounting}.

\subsection{The resolved conifold case}
The web diagram is
$$\xy
(-5, 0); (-5, 5), **@{-};
(-5, 0); (-10, 0), **@{-};
(-5, 0); (0, -5), **@{-};
(0, -5); (5, -5), **@{-};
(0, -5); (0, -10), **@{-}; (-1, -1)*+{\mu};
\endxy$$
\begin{center}
{\em Figure 1}
\end{center}
\vspace{.1in}
The (undeformed) partition function is given by:
\begin{eqnarray*}
Z_{(0)} & = & \sum_{\mu} \cW_{\mu^t}(q) (-1)^{|\mu|} Q^{|\mu|} \cW_{\mu}(q).
\end{eqnarray*}
We take the deformed partition function to be:
\begin{eqnarray}
\cZ_{(0)} & = & \sum_{\mu} e^{-|\mu|t} W_{\mu^t}(q, e^{-t}) (-1)^{|\mu|} Q^{|\mu|} W_{\mu}(q, e^{-t}).
\end{eqnarray}
By (\ref{eqn:WMut2}) it can be also written as:
\begin{eqnarray}
\cZ_{(0)} & = & \sum_{\mu} e^{-|\mu|t}  Q^{|\mu|} W_{\mu}(q, e^{-t})W_{\mu^t}(q, e^{t}).
\end{eqnarray}
By (\ref{def:WMu}),
one easily gets:
\begin{eqnarray} \label{eqn:ZOrbifold}
&& \cZ_{(0)} = \sum_{\mu} \prod_{(i, j) \in \mu} \frac{(1 - e^{-t} q^{j-i})(1- e^{-t} q^{i-j})}
{(1 - q^{h(i, j)})(1 - q^{-h(i,j)})}.
\end{eqnarray}

We will prove two results for such deformed partition functions.
First of all,
we will identify it with certain partition functions on the Hilbert schemes.
Secondly,
we will find infinite product expression for it so that one can extract the {\em deformed
Gopakumar-Vafa invariants}.

By comparing (\ref{eqn:ZOrbifold}) and (\ref{eqn:ChiVV-}),
one gets the following result which generalizes the undeformed case
considered in \cite{Hol-Iqb-Vaf, Li-Liu-Zho}:

\begin{theorem}
One has the following identification of deformed partition functions:
\begin{eqnarray}
&& \cZ_{(0)} = \sum_{n=0}^{\infty} Q^n
\chi((\bC^2)^{[n]}, \Lambda_{-e^{-t}}(V_n) \otimes \Lambda_{-e^{-t}}(V_n^*))(q, q^{-1}).
\label{eqn:Z0}
\end{eqnarray}
\end{theorem}

The second result is a straightforward consequence of (\ref{eqn:L00}) and the observation
\begin{eqnarray*}
\cZ_{(0)} & = &  \cL_{(0)(0)}(Q).
\end{eqnarray*}
It can be stated as follows.

\begin{theorem}
The deformed partition function for the resolved conifold can be rewritten as follows.
\begin{eqnarray}
\cZ_{(0)} & = & \exp (-\sum_{n=1}^{\infty} \frac{(1 - e^{-nt})^2}{(q^{n/2} - q^{-n/2})^2}Q^n)
\label{eqn:Z02} \\
& = & \prod_{n=1}^{\infty} \frac{(1 - Qq^n)^n(1 - Qe^{-t}q^n)^n}{(1 - Qe^{-t/2}q^n)^{2n}}.
\end{eqnarray}
\end{theorem}

Besides the natural $T^2$-action on $\bC^2$,
consider an extra circle action on $V_1$ and $V_1^*$ by scalar multiplications along the fibers.
Then we have
$$\chi(\bC^2, \Lambda_{-1}(V_1) \otimes \Lambda_{-1}(V_1^*))(t_1, t_2, e^{-t})
= \frac{(1-e^{-t})^2}{(1-t_1^{-1})(1-t_2^{-1})}.$$
Then we have
\begin{eqnarray}
&& \cZ_{(0)} = \exp \sum_{n=1}^{\infty} \frac{1}{n}
\chi(\bC^2, \Lambda_{-1}(V_1) \otimes \Lambda_{-1}(V_1^*))(t_1^n, t_2^n, e^{-nt})|_{t_1 = q, t_2= q^{-1}}.
\end{eqnarray}
From this one expects an interpretation in terms of symmetric products
generalizing the nondeformed case in \cite{Li-Liu-Zho}.

\subsection{A related case}
A related case is given by the following diagram,
where the two horizontal edges are glued with each other.
$$\xy
(-5, 0); (-5, 5), **@{-};
(-5, 0); (-10, 0), **@{-};
(-5, 0); (0, -5), **@{-};
(0, -5); (5, -5), **@{-};
(0, -5); (0, -10), **@{-};
(-7.5, 0)*+{|}; (2.5, -5)*+{|};
\endxy$$
\begin{center}
{\em Figure 2}
\end{center}
\vspace{.1in}
The (undeformed) partition function is given by
\begin{eqnarray}
&& Z(Q, Q_m, q) = \sum_{\mu, \nu } (-Q)^{|\mu|}
(-Q_m)^{|\nu|} \cW_{\mu\nu}(q) \cW_{\nu^t\mu^t}(q).
\end{eqnarray}
We take the deformed partition function to be
\begin{eqnarray} \label{drf:Z5D}
&& \tilde{\cZ}(Q, Q_m, q, e^{-t}) = \sum_{\mu, \nu } (-Q)^{|\mu|}
(-Q_m)^{|\nu|}  W_{\mu\nu}(q, e^{-t})W_{\nu^t\mu^t}(q, e^{-t}).
\end{eqnarray}
Now note that
\begin{eqnarray}
&& \tilde{\cZ}(Q, Q_m, q, e^{-t})
= \sum_{\mu} (-Q)^{|\mu|} \cL_{\mu(\mu)^t}(Q_m).
\end{eqnarray}
Hence by (\ref{eqn:LMuMu})
\begin{equation}
\begin{split}
& \frac{\tilde{\cZ}(Q, Q_m, q, e^{-t})}{\cZ_{(0)}(Q_m, q, e^{-t})} \\
= & \sum_{\mu} Q^{|\mu|} \prod_{e \in \mu}
\frac{(1 - e^{-t} q^{c(e)})(1 - e^{-t} q^{-c(e)})}{(1- Q_me^{-t}q^{c(e)})(1-Q_me^{-t}q^{-c(e)})} \\
& \cdot \prod_{e \in \mu} \frac{(1 - Q_mq^{h(e)})(1 - Q_mq^{-h(e)})}{(1 - q^{h(e)})(1 - q^{-h(e)})}.
\end{split}
\end{equation}

We now identify this as a partition function on the Hilbert schemes by comparing
with (\ref{eqn:ChiWq}),
generalizing the results in the deformed case \cite{Hol-Iqb-Vaf, Li-Liu-Zho}:

\begin{theorem}
One has the following identification of deformed partition functions:
\begin{eqnarray}
&& \frac{\tilde{\cZ}(Q, Q_m)}{\cZ_{(0)}(Q_m)} = \sum_{n=0}^{\infty} Q^n
\chi((\bC^2)^{[n]}, \Lambda_{-Q_m}(\bV_n) \otimes W_{e^{-t}, Q_m}(\bV_n^*))(q, q^{-1}).
\end{eqnarray}
\end{theorem}

We do not know a rigorous method to establish an infinite product expression for $\tilde{\cZ}(Q, Q_m)$ at present.
For the unddeformed partition function,
a physical calcultion based on duality with Chern-Simons thoery was carried out in \cite{Hol-Iqb-Vaf},
and a mathematical proof based on symmetric products was given in \cite{Li-Liu-Zho}.
Here we first present another mathematical proof for the undeformed case using Schur calculus.

\begin{theorem} \cite{Hol-Iqb-Vaf, Li-Liu-Zho}
For the undeformed partition function,
one has
\begin{eqnarray}
&& Z(Q, Q_m, q) = \prod_{k=0}^{\infty} \prod_{l=1}^{\infty}
\frac{(1- Q^{k+1}Q_m^kq^l)^l(1- Q^kQ_m^{k+1}q^l)^l}{(1- Q^{k+1}Q_m^{k+1}q^l)^{2l}}.
\end{eqnarray}
\end{theorem}

\begin{proof}
We use the following identity \cite{ZhoConjecture}:
\begin{eqnarray}
&& W_{\mu\nu}(q) = (-1)^{|\mu|+\nu|} q^{(\kappa_{\mu}+\kappa_{\nu})/2}
\sum_{\eta} s_{\mu/\eta}(q^{-\rho})s_{\nu/\rho}(q^{-\rho}).
\end{eqnarray}
and (\ref{eqn:Schur4SumT}):
\begin{eqnarray*}
&& Z(Q, Q_m, q) = \sum_{\mu, \nu } (-Q)^{|\mu|}
(-Q_m)^{|\nu|}  W_{\mu\nu}(q)W_{\nu^t\mu^t}(q) \\
& = &  \sum_{\mu, \nu, \xi, \eta } (-Q)^{|\mu|}
(-Q_m)^{|\nu|}  s_{\mu/\xi}(q^{-\rho}) s_{\nu/\xi}(q^{-\rho})
s_{\nu^t/\eta}(q^{-\rho})s_{\mu^t/\eta}(q^{-\rho}) \\
& = &  \sum_{\mu, \nu, \xi, \eta } Q^{|\mu|}
Q_m^{|\nu|}  s_{\mu/\xi}(q^{-\rho}) s_{\nu/\xi}(q^{-\rho})
s_{\nu/\eta}(q^{\rho})s_{\mu/\eta}(q^{\rho}) \\
& = & \exp \sum_{n=1}^{\infty}
\frac{1}{n(1 - (QQ_m)^n)}
(Q^np_n(q^{-\rho})p_n(q^{\rho})
+ Q_m^np_n(q^{-\rho})p_n(q^{\rho}) \\
&& + (QQ_m)^np_n(q^{-\rho})p_n(q^{-\rho})
+ (QQ_m)^n p_n(q^{\rho})p_n(q^{\rho}) ).
\end{eqnarray*}
By (\ref{eqn:pn}) we then have:
\begin{eqnarray*}
&& Z(Q, Q_m, q)\\
& = & \exp \sum_{n=1}^{\infty}
\frac{1}{n(1 - (QQ_m)^n)} ( - \frac{Q^nq^n}{(1- q^n)^2} - \frac{Q_m^nq^n}{(1-q^n)^2}
+ \frac{2 (QQ_m)^nq^n}{(1-q^n)^2}) \\
& = & \prod_{k=0}^{\infty} \prod_{l=1}^{\infty} \frac{(1- Q^{k+1}Q_m^kq^l)^l(1- Q^kQ_m^{k+1}q^l)^l}
{(1- Q^{k+1}Q_m^{k+1}q^l)^{2l}}.
\end{eqnarray*}
\end{proof}

It turns out that this proof can be generalized to another deformed partition function:
\begin{eqnarray} \label{def:Z5D2}
&& \cZ(Q, Q_m; q, e^{-t}) = \sum_{\nu, \eta} (-Q_1)^{|\nu|} (-Q_m)^{|\eta|}
W_{\eta(0)\nu}(q, e^{-t})W_{\eta^t(0)\nu^t}(q, e^{-t}).
\end{eqnarray}

\begin{theorem}
We have
\begin{equation}
\begin{split}
&  \cZ(Q, Q_m; q, e^{-t}) \\
= & \prod_{k=0}^{\infty}\left( \prod_{l=1}^{\infty}
\left(\frac{(1- Q^{k+1}Q_m^k q^l)(1- e^{-2t}Q^{k+1}Q_m^k q^l)}
{(1- e^{-t}Q^{k+1}Q_m^kq^l)^2}\right)^l  \right.\\
& \cdot \left(\frac{(1- Q^{k}Q_m^{k+1} q^l)(1- e^{-2t}Q^{k}Q_m^{k+1} q^l)}
{(1- e^{-t}Q^{k}Q_m^{k+1})^2}\right)^l \\
& \left. \cdot \left(\frac{(1- e^{-t}Q^{k+1}Q_m^{k+1}q^l)^2}
{(1- Q^{k+1}Q_m^{k+1} q^l)(1- e^{-2t}Q^{k+1}Q_m^{k+1} q^l)}\right)^{2l} \right).
\end{split}
\end{equation}
\end{theorem}

\begin{proof}
By (\ref{eqn:Schur4SumT}) we have:
\begin{eqnarray*}
&&  \cZ(Q, Q_m; q, e^{-t}) \\
& = & \sum_{\nu, \eta} (-Q)^{|\nu|} (-Q_m)^{|\eta|} \\
&& \cdot \sum_{\xi^1} q^{(|\eta| - |\xi^1|)/2} s_{\eta/\xi^1}(E(u; q, e^{-t}))
q^{(|\nu| - |\xi^1|)/2} s_{\nu^t/\xi^1}(E(u; q, e^{-t})) \\
&& \cdot \sum_{\xi^2}
q^{(|\eta^t| - |\xi^2|)/2} s_{\eta^t/\xi^2}(E(u; q, e^{-t}))
q^{(|\nu| - |\xi^2|)/2} s_{\nu/\xi^2}(E(u; q, e^{-t})) \\
& = & \exp \sum_{n=1}^{\infty} \frac{q^n}{n(1- (Q_1Q_m)^n)} \\
&& (-Q^n - Q_m^n + 2 (QQ_m)^n ) p_n(E(u; q, e^{-t}))p_n(E(u; q, e^{-t})) \\
& = & \exp \sum_{n=1}^{\infty} \frac{q^n}{n(1- (Q_1Q_m)^n)} \frac{(1- e^{-nt})^2}{(q^n-1)^2}
(-Q^n - Q_m^n + 2 (QQ_m)^n )
\end{eqnarray*}
It is then an easy exercise to show that:
\begin{eqnarray*}
\cZ(Q, Q_m; q, e^{-t})
& = & \prod_{k=0}^{\infty}\left( \prod_{l=1}^{\infty}
\left(\frac{(1- Q^{k+1}Q_m^k q^l)(1- e^{-2t}Q^{k+1}Q_m^k q^l)}
{(1- e^{-t}Q^{k+1}Q_m^kq^l)^2}\right)^l  \right.\\
&& \cdot \left(\frac{(1- Q^{k}Q_m^{k+1} q^l)(1- e^{-2t}Q^{k}Q_m^{k+1} q^l)}
{(1- e^{-t}Q^{k}Q_m^{k+1})^2}\right)^l \\
&& \left. \cdot \left(\frac{(1- e^{-t}Q^{k+1}Q_m^{k+1}q^l)^2}
{(1- Q^{k+1}Q_m^{k+1} q^l)(1- e^{-2t}Q^{k+1}Q_m^{k+1} q^l)}\right)^{2l} \right).
\end{eqnarray*}
\end{proof}

\subsection{Another related case}

Consider the local Calabi-Yau geometry with the following web diagram:
$$
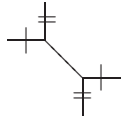
\xy
(-5, 0); (-5, 5), **@{-};
(-5, 0); (-10, 0), **@{-};
(-5, 0); (0, -5), **@{-};
(0, -5); (5, -5), **@{-};
(0, -5); (0, -10), **@{-};
(-7.5, 0)*+{|}; (2.5, -5)*+{|}; (-4.7, 2.5)*+{=}; (0, -7.5)*+{=};
\endxy$$
where we glue together the horizontal edge with other and also the vertex edge with each other.
The undeformed theory of this case has been discussed in Section 5.1.1 and Section 6.2 in \cite{Hol-Iqb-Vaf}.
The partition function by topological vertex method is
\begin{eqnarray}
&& Z(Q, Q_m, Q_1; q) = \sum_{\mu, \nu, \eta} (-Q)^{|\mu|} (-Q_1)^{|\nu|} (-Q_m)^{|\eta|}
\cW_{\mu\nu\eta}(q)\cW_{\mu^t\nu^t\eta^t}(q).
\end{eqnarray}
We take the deformed partition function to be:
\begin{equation} \label{def:Z6D}
\begin{split}
& \cZ(Q, Q_m, Q_1; q, e^{-t}) \\
= & \sum_{\mu, \nu, \eta}
(-Q)^{|\mu|} (-Q_1)^{|\nu|} (-Q_m)^{|\eta|}
W_{\eta\mu\nu}(q, e^{-t})W_{\eta^t\mu^t\nu^t}(q, e^{-t}).
\end{split}
\end{equation}

\begin{theorem}
We have
\begin{equation} \label{eqn:Z6D}
\begin{split}
&  \frac{\cZ(Q, Q_m, Q_1; q, e^{-t})}{\tilde{\cZ}(Q_1, Q_m; q, e^{-t})} \\
= & \sum_{\mu} Q^{|\mu|} \prod_{(i, j) \in \mu}
\frac{(1 - e^{-t} q^{j-i})(1 - e^{-t} q^{-(j-i)})}{(1- q^{h(i,j)})(1-q^{-h(i,j)})} \\
& \prod_{k=1}^{\infty}
\frac{(1 - Q_m Q_{\rho}^{k-1} q^{h(i,j)})(1 - Q_mQ_{\rho}^{k-1} q^{-h(i, j)})}
{(1 - e^{-t} Q_m Q_{\rho}^{k-1} q^{j - i})(1 - e^{-t} Q_mQ_{\rho}^{k-1} q^{-(j - i)})} \\
& \frac{(1 - Q_1 Q_{\rho}^{k-1} q^{h(i,j)})(1 - Q_1Q_{\rho}^{k-1} q^{-h(i, j)})}
{(1 - e^{-t} Q_1Q_{\rho}^{k-1}q^{j - i})(1 - e^{-t} Q_1Q_{\rho}^{k-1} q^{-(j - i)})} \\
& \frac{(1 - Q_{\rho}^ke^{-t} q^{j - i})^2(1 - Q_{\rho}^k e^{-t} q^{-(j - i)})^2}
{(1 - Q_{\rho}^k q^{h(i,j)})^2(1 - Q_{\rho}^k q^{-h(i, j)})^2}.
\end{split}
\end{equation}
\end{theorem}

\begin{proof}
First of all,
by the definition of the deformed topological vertex we have:
\begin{eqnarray*}
&& \cZ(Q, Q_m, Q_1; q, e^{-t}) \\
& = & \sum_{\mu, \nu, \eta} (-Q)^{|\mu|}(-Q_1)^{|\nu|} (-Q_m)^{|\eta|}
e^{-|\mu|t/2}W_{\mu}(q, e^{-t}) e^{-|\mu^t|t/2} W_{\mu^t}(q, e^{-t}) \\
&& \cdot \sum_{\xi^1}
q^{(|\eta| - |\xi^1|)/2} s_{\eta/\xi^1}(E^{\mu^t}(u; q, e^{-t}))
q^{(|\nu| - |\xi^1|)/2} s_{\nu^t/\xi^1}(E^{\mu}(u; q, e^{-t})) \\
&& \cdot \sum_{\xi^2}
q^{(|\eta^t| - |\xi^2|)/2} s_{\eta^t/\xi^2}(E^{\mu}(u; q, e^{-t}))
q^{(|\nu| - |\xi^2|)/2} s_{\nu/\xi^2}(E^{\mu^t}(u; q, e^{-t})).
\end{eqnarray*}
Hence by (\ref{eqn:Schur4SumT})  we have:
\begin{eqnarray*}
&& \cZ(Q, Q_m, Q_1; q, e^{-t}) \\
& = & e^{-|\mu|t}W_{\mu}(q, e^{-t}) W_{\mu^t}(q, e^{-t})
\exp \sum_{n=1}^{\infty} \frac{q^n}{n(1- (Q_1Q_m)^n)} \\
&& \cdot (-Q_1^n - Q_m^n + 2 (Q_1Q_m)^n )
p_n(E^{\mu}(u; q, e^{-t}))p_n(E^{\mu^t}(u; q, e^{-t})).
\end{eqnarray*}
It follows that
\begin{eqnarray*}
&& \frac{\cZ(Q, Q_m, Q_1; q, e^{-t})}{\cZ(Q_1, Q_m; q, e^{-t})} \\
& = & e^{-|\mu|t}W_{\mu}(q, e^{-t}) W_{\mu^t}(q, e^{-t})
\exp \sum_{n=1}^{\infty} \frac{q^n}{n(1- (Q_1Q_m)^n)} \\
&& \cdot (-Q_1^n - Q_m^n + 2 (Q_1Q_m)^n ) F_{\mu, \mu^t}(q^n, e^{-nt}).
\end{eqnarray*}
Hence (\ref{eqn:Z6D}) follows from (\ref{eqn:FProd3}).
\end{proof}

\begin{theorem}
We have the following identity:
\begin{equation} \label{eqn:6DIden}
\begin{split}
& \frac{\cZ(Q, Q_m, Q_1; q, e^{-t})}{\cZ(Q_1, Q_m; q, e^{-t})} \\
= & \sum_{n \geq 0} Q^n \chi((\bC^2)^{[n]}, \Ell(T(\bC^2)^{[n]}, \bV_n)(Q_m, Q_{\rho}, e^{-t}))(q, q^{-1}),
\end{split}
\end{equation}
where $Q_{\rho} = Q_1Q_m$.
\end{theorem}

\begin{proof}
we rewrite the right-hand side of (\ref{eqn:Z6D}) as follows.
\begin{eqnarray*}
&& \sum_{\mu}Q^{|\mu|} \prod_{(i, j)\in \mu}
\prod_{k=1}^\infty
\left( \frac{(1-Q_\rho^{k-1}Q_mq^{h(i, j)})(1-Q_\rho^{k-1}Q_mq^{-h(i, j)})}
{(1-Q_\rho^{k-1}q^{h(i,j)})(1-Q_\rho^{k-1}q^{-h(i, j)})} \right. \\
&& \cdot \frac{(1-Q_\rho^kQ_m^{-1} q^{h(i, j)})(1-Q_\rho^kQ_m^{-1}q^{-h(i, j)})}
{(1-Q_\rho^kq^{h(i,j)})(1-Q_\rho^kq^{-h(i, j)})} \\
&& \cdot \frac{(1-e^{-t} Q_\rho^{k-1}q^{c(i,j)})(1-e^{-t}Q_\rho^{k-1}q^{-c(i, j)})}
{(1-e^{-t} Q_\rho^{k-1}Q_mq^{c(i, j)})(1-e^{-t}Q_\rho^{k-1}Q_mq^{-c(i, j)})} \\
&& \left. \cdot \frac{(1-e^{-t}Q_\rho^kq^{c(i,j)})(1-e^{-t}Q_\rho^kq^{-c(i, j)})}
{(1-e^{-t}Q_\rho^kQ_m^{-1} q^{c(i, j)})(1-e^{-t}Q_\rho^kQ_m^{-1}q^{-c(i, j)})} \right),
\end{eqnarray*}
Then (\ref{eqn:6DIden}) follows from (\ref{eqn:ChiEllq}).
\end{proof}

\subsection{The $SU(2)$ case}
The web diagrams are:
$$
\xy
(10, 0); (5, 5), **@{-}; (15, 2)*+{B}; (45, 2)*+{B}; (84, 2)*+{B};
(10,0); (20,0), **@{-}; (25, 5); **@{-}; (8, -5)*+{F}; (78, -5)*+{F}; (38,-5)*+{F};
(10,-10); (20,-10), **@{-}; (25, -15), **@{-}; (15, -12)*+{B}; (50, -12)*+{B+F}; (95, -12)*+{B+2F};
(10,0); (10,-10), **@{-}; (5, -15), **@{-};
(20,0); (20,-10), **@{-};  (15, -5)*+{\bF_0}; (22, -5)*+{F}; (57, -5)*+{F}; (107, -5)*+{F};
(40, 0); (35, 5), **@{-}; (45, 2)*+{B};
(40, 0); (50, 0), **@{-}; (60, -10),  **@{-}; (40, -10),  **@{-};
(40, 0),  **@{-}; (47, -5)*+{\bF_1}; (38,-5)*+{F};
(40, -10); (35, -15),  **@{-};
(50, 0); (50, 5),  **@{-};
(60, -10); (70, -15),  **@{-};
(75, 5); (80, 0), **@{-}; (90, 0), **@{-}; (110, -10),  **@{-}; (80, -10),  **@{-};
(80, 0),  **@{-}; (90, -5)*+{\bF_2};
(80, -10); (75, -15),  **@{-};
(90, 0); (85, 5),  **@{-};
(110, -10); (125, -15),  **@{-};
\endxy$$
\begin{center}
{\em Figure 3}
\end{center}
\vspace{.1in}
where $B, F \in H_2(\bF_m, \bZ)$ are the homological classes of the base and the fiber respectively,
and we have
\begin{align*}
F^2 & = 0, & B^2 & = -m, & BF & = 1,
\end{align*}
hence
$$(B+mF)^2 = m.$$
The undeformed partition functions are calculated in physics literature in
\cite{Aga-Mar-Vaf} and formulated as a Feynman rule in \cite{Iqb}.
The mathematical proof can be found in \cite{ZhoChemistry}.
The self-intersection numbers of the divisors in the surfaces represented by the internal
edges play a role in the formula:
\begin{eqnarray*}
Z_{\bF_m}(Q_B, Q_F) & = &
\sum_{\mu^{1, 2}, \nu^{1, 2}} \cW_{\mu^1\nu^1}(q)Q_F^{|\nu^1|} \cW_{\nu^1(\mu^2)^t}(q)
(-1)^{m|\mu^2|}q^{m\kappa_{(\mu^2)^t}} \\
&& \cdot (Q_BQ_F^m)^{|\mu^2|} \cW_{(\mu^2)^t\nu^2}(q)Q_F^{|\nu^2|} \cW_{\nu^2\mu^1}(q)
(-1)^{-m|\mu^1|}q^{-m\kappa_{\mu^1}} \\
& = &  \sum_{\mu^1, \mu^2} Q_B^{|\mu^1|+|\mu^2|}
[(-1)^{|\mu^1|+|\mu^2|}Q_F^{|\mu^2|}q^{-\frac{1}{2}(\kappa_{\mu^1}+\kappa_{\mu^2})}]^m
K_{\mu^1, (\mu^2)^t}(Q_F)^2.
\end{eqnarray*}
Here the locations of the partitions are indicated below:
$$
\xy
(75, 5); (80, 0), **@{-}; (90, 0), **@{-}; (110, -10),  **@{-}; (80, -10),  **@{-};
(80, 0),  **@{-}; (90, -5)*+{\bF_m};
(80, -10); (75, -15),  **@{-};
(90, 0); (85, 5),  **@{-};
(110, -10); (125, -15),  **@{-};
(84, 2)*+{\mu^1}; (78, -5)*+{\nu^1}; (93, -13)*+{(\mu^2)^t}; (101, -2)*+{\nu^2};
\endxy$$
\begin{center}
{\em Figure 4}
\end{center}
\vspace{.1in}
The above expression for $Z_{\bF_m}$ suggests the following partition function:
\begin{equation} \label{eqn:ZFm}
\begin{split}
& \tilde{\cZ}_{\bF_m}(q, Q_B, Q_F, e^{-t}) \\
= &
\sum_{\mu^{1, 2}, \nu^{1, 2}}
e^{-(2|\mu^1|+2|\mu^2|+|\nu^1|+|\nu^2|)t} \\
& \cdot W_{\mu^1\nu^1}(q, e^{-t})Q_F^{|\nu^1|} W_{\nu^1(\mu^2)^t}(q, e^{-t})
(-1)^{m|\mu^2|}q^{m\kappa_{(\mu^2)^t}} \\
& \cdot (Q_BQ_F^m)^{|\mu^2|} W_{(\mu^2)^t\nu^2}(q, e^{t})Q_F^{|\nu^2|} W_{\nu^2\mu^1}(q, e^{t})
(-1)^{-m|\mu^1|}q^{-m\kappa_{\mu^1}}.
\end{split}
\end{equation}
This is not a deformation of $Z_{\bF_m}$.
However it has a similar property as $Z_{\bF_m}$ as follows.

\begin{theorem}
We have
\begin{equation} \label{eqn:NekN=2}
\begin{split}
& \frac{\cZ_{\bF_m}(Q_B, Q_F; q, e^{-t})}{\cK_{(0)(0)}(Q_F; q, e^{-t})\cK_{(0)(0)}(Q_F; q, e^{t})} \\
= & \sum_{k \geq 0} Q_B^k \chi(M(2, k), \Lambda_{-e^{-t}}(\bV_k \otimes \bW_k^*) \otimes
\Lambda_{-e^{-t}}(\bV_k^* \otimes \bW_k) \otimes E_{N, k}^m)(e_1, e_2, q, q^{-1})
\end{split} \end{equation}
where $Q_F = e_1e_2^{-1}$, $e_1 = -1$.
\end{theorem}

\begin{proof}
By (\ref{eqn:ZFm}) we have
\begin{eqnarray*}
&& \frac{\cZ_{\bF_m}(q, Q_B, Q_F, e^{-t})}{\cK_{(0)(0)}(Q_F; q, e^{-t})\cK_{(0)(0)}(Q_F; q, e^t)} \\
& = &  \sum_{\mu^1, \mu^2} Q_B^{|\mu^1|+|\mu^2|}
[(-1)^{|\mu^1|+|\mu^2|}Q_F^{|\mu^2|}q^{-\frac{1}{2}(\kappa_{\mu^1}+\kappa_{\mu^2})}]^m\\
&& \cdot \frac{\cK_{\mu^1, (\mu^2)^t}(Q_F; q, e^{-t})\cK_{\mu^1, (\mu^2)^t}(Q_F; q, e^{t})}
{\cK_{(0)(0)}(Q_F; q, e^{-t})\cK_{(0)(0)}(Q_F; q, e^{t})}.
\end{eqnarray*}
Now we can apply Theorem \ref{thm:K} and (\ref{eqn:Rank>1}).
\end{proof}

By the second product expression in Theorem \ref{thm:K} one can also express
$$\frac{\cZ_{\bF_m}(q, Q_B, Q_F, e^{-t})}{\cK_{(0)(0)}(Q_F; q, e^{-t})\cK_{(0)(0)}(Q_F; q, e^t)}$$
as a sum of infinite products.
This expression generalizes Nekrasov's partition function.

We propose the deformed partition function to be:
\begin{eqnarray*}
&& \cZ_{\bF_m}(Q_B, Q_F; q, e^{-t}) \\
& = & \sum_{\mu^{1, 2}, \nu^{1, 2}} (-1)^{m(|\mu^1|+|\mu^2|)} q^{-m\kappa_{\mu^1}}
q^{m\kappa_{\mu^2}/2}
Q_B^{|\mu^1|}(Q_BQ_F^m)^{|\mu^2|} \\
&& \cdot W_{\mu^1\nu^1}(q, e^{-t})  Q_F^{|\nu^1|} W_{\nu^1\mu^2}(q, e^{-t})
\cdot W_{\mu^2\nu^2}(q, e^{-t}) Q_F^{|\nu^2|} W_{\nu^2\mu^1}(q, e^{-t}).
\end{eqnarray*}
We have checked in some low degree cases that this to satisfy the following deformed version of Gopakumar-Vafa integrality:
\begin{eqnarray}
F & = & \sum_{\Sigma \in H_2(X)-\{0\}} \sum_{g \geq 0} \sum_{k \geq 1} \frac{1}{k}
(-1)^{g-1}n^g_{\Sigma}(e^{-kt}) (q^{\frac{k}{2}} - q^{-\frac{k}{2}})^{2g-2}
Q^{k\Sigma},
\end{eqnarray}
where $n^g_{\Sigma}(e^{-t})$ is a polynomial in $e^{-t}$ with integral coefficients.
For example, for $\bF_0$,
we have
\begin{eqnarray*}
&& n^g_{B}(e^{-t}) = n^g_{F} = -2 (1 - e^{-t})^2 \delta_{g, 0}, \\
&& n^g_{B^2} = n^g_{F^2} = 0, \;\;\;\; n^g_{B+F} = -4 (1 - e^{-t})^3, \\
&& n^g_{B^3} = n^g_{F^3} = 0, \;\;\;\; n^g_{2B+F} = n^g_{B+2F} = - (6 - 2e^{-t})(1 - e^{-t})^3\delta_{g, 0}.
\end{eqnarray*}
For $\bF_1$ we have
\begin{eqnarray*}
&& n^g_{B} = (1-e^{-t})^2 \delta_{g, 0}, \;\;\;\;
n^g_F = -2 (1 - e^{-t})^2 \delta_{g, 0}, \\
&& n^g_{2B} = n^g_{2F} = 0, \;\;\;\;
n^g_{B+F} = (1 - e^{-t})^3(3 + e^{-t}) \delta_{g, 0}, \\
&& n^0_{3B} = -2e^{-3t}(1- e^{-t})^3 \delta_{g, 0}, \;\;\;\;
n^1_{B^3} = - e^{-4t}, \\
&& n^g_{2B+F} = (1-e^{-t})^2(2e^{-2t} + 2e^{-3t}) \delta_{g, 0}, \\
&& n^g_{B+2F} = (1 - e^{-t})^2(5-6e^{-t}-+e^{-2t})\delta_{g, 0}.
\end{eqnarray*}
For $\bF_2$ we have
\begin{eqnarray*}
&& n^g_B = -(1-e^{-t})^2\delta_{g, 0}, \;\;\;\; n^g_F = -2(1 - e^{-t})^2 \delta_{g, 0}, \\
&& n^0_{2B} = 2e^{-t}(1-e^{-t})^3, \;\;\;\; n^0_{2B} = e^{-t}(1-e^{-t})^2, \\
&& n^g_{B+F} = - 2(1 - e^{-t})^3 \delta_{g, 0}, \;\;\;\; n^g_{2F} = 0.
\end{eqnarray*}
In general,
one can consider the deformed topological string partition functions for toric local Calabi-Yau
$3$-folds in the same fashion and study the deformed Gopakumar-Vafa invariants.
For example,
for $\bP^2$ we have
\begin{eqnarray*}
&& n^g_1(e^{-t}) = 3 (1 - e^{-t})^2 \delta_{g, 0}, \\
&& n^g_2(e^{-t}) = -3 (1 - e^{-t})^2(2 - e^{-t}) \delta_{g, 0}, \\
&& n^0_3(e^{-t}) = (27-15e^{-t})(1- e^{-t})^3, \;\;\; n^1_3 = (-10+6e^{-t})(1 - e^{-t})^2, \;\;\;\;
n^g_3 = 0 (g > 1), \\
&& n^0_4 = -(1-e^{-t})^2(192-474e^{-t}+390e^{-2t}-114e^{-3t} + 6e^{-4t}), \\
&& n^1_4 = (1-e^{-t})^2(231-402e^{-t}+201e^{-2t} - 24e^{-3t}), \\
&& n^2_4 = -(1-e^{-t})^2(102-120e^{-t}+30e^{-2t}), \\
&& n^3_4 = (1 - e^{-t})^2 (15- 12 e^{-t}), \\
&& n^g_4 = 0, \;\;\;\; g > 3.
\end{eqnarray*}
We notice the following positivity of $n^g_d$ for these example.
Define
$$P^g_d(x) = (-1)^{g+d-1}n^g_d(-x).$$
Then the coefficients of $P^g_d(x)$ are nonnegative.
We conjecture it is true in general.

Here we want to offer some speculations on the geometric meaning of extra variable in the partition function.
For a toric Fano surface $X$,
the mathematical definition of its Gromov-Witten invariants are defined as follows:
$$\int_{[\Mbar_{g, 0}(X, d)]^{\virt}} e(K_X^{g,d}),$$
where $d  \in H_2(X, \bZ)$, $\Mbar_{g, 0}(X, d)$ is the moduli space
of genus $g$ stable maps of degree $d$ to $X$,
$K_X^{g,d}$ is the vector bundle on $\Mbar_{g, 0}(X, d)$ whose fiber at a map $f: \Sigma \to X$
is given by $H^1(\Sigma, f^*K_X)$.
Since $\Mbar_{g, 0}(X, d)$ is compact,
one can do the calculations in the equivariant setting:
$$\int_{[\Mbar_{g, 0}(X, d)]_{T^2}^{\virt}} e_{T^2}(K_X^{g,d}),$$
and the result is a constant.
Here we use the $T^2$-action that defines the toric structure.
However,
it is tempting also to do the calculation equivariantly directly on $K_X$.
For this purpose consider an extra circle action on the fiber of the $K_X$ by
multiplication and consider
$$\int_{[\Mbar_{g, 0}(K_X, d)]_{T^3}^{\virt}} 1.$$
One can take the result obtained by formally applying the localization as definition.
This might explain the extra parameter in the deformed partition functions.

\subsection{The $SU(N)$ ($N>2$) cases}

The undeformed topological string partition function is given by the topological
vertex as follows \cite{Iqb-Kas2}:
\begin{equation} \label{eqn:N>2}
\begin{split}
Z^{(m)}_{A_{N-1}} = & \sum_{\mu^1, \dots, \mu^N}
\prod_{i=1}^N (Q_{b_i}^{|\mu^i|} (-1)^{(N+m-2i)|\mu^i|} q^{(N+2m-2i+1)\kappa_{\mu^i}}) \\
& \cdot \sum_{\nu^{1, \dots, N-1}}
\prod_{k=1}^N (\cW_{\nu^{k-1}\mu^k(\nu^k)^t}(q) q^{\kappa_{\nu^k/2}} Q_k^{|\nu^k|}) \\
& \cdot \sum_{\eta^{1, \dots, N-1}}
\prod_{k=1}^N (\cW_{\eta^{k-1}(\mu^{N+1-k})^t(\eta^k)^t}(q) q^{\kappa_{\eta^k/2}} Q_{N+1-k}^{|\eta^k|}),
\end{split}
\end{equation}
where $\nu^0 = \nu^N = \eta^0=\eta^N = (0)$, $Q_N=1$.
This suggests the following deformed partition function:
\begin{equation} \label{eqn:N>2Def}
\begin{split}
\cZ^{(m)}_{A_{N-1}} = & \sum_{\mu^1, \dots, \mu^N}
\prod_{i=1}^N (Q_{b_i}^{|\mu^i|} (-1)^{(N+m-2i)|\mu^i|} q^{(N+2m-2i+1)\kappa_{\mu^i}}) \\
& \cdot \sum_{\nu^{1, \dots, N-1}}
\prod_{k=1}^N (W_{\nu^{k-1}\mu^k(\nu^k)^t}(q, e^{-t}) q^{\kappa_{\nu^k/2}} Q_k^{|\nu^k|}) \\
& \cdot \sum_{\eta^{1, \dots, N-1}}
\prod_{k=1}^N (W_{\eta^{k-1}(\mu^{N+1-k})^t(\eta^k)^t}(q, e^{-t}) q^{\kappa_{\eta^k/2}} Q_{N+1-k}^{|\eta^k|}),
\end{split}
\end{equation}
where $\nu^0 = \nu^N = \eta^0=\eta^N = (0)$, $Q_N=1$.
We expect this to satisfy the deformed Gopakumar-Vafa integrality.

On the other hand, one can use the following symmetry of the topological vertex (see e.g. \cite{ZhoCounting}):
$$\cW_{(\mu^1)^t(\mu^2)^t(\mu^3)^t}(q) = q^{-(\kappa_{\mu^1}+\kappa_{\mu^2}+\kappa_{\mu^3})/2}
\cW_{\mu^3,\mu^2, \mu^1}(q)$$
to rewrite the last line in (\ref{eqn:N>2}).
This yields the following expression \cite{Iqb-Kas2}:
\begin{equation} \label{eqn:N>2New}
\begin{split}
Z^{(m)}_{A_{N-1}} = & \sum_{\mu^1, \dots, \mu^N}
\prod_{i=1}^N (Q_{b_i}^{|\mu^i|} (-1)^{(N+m)|\mu^i|} q^{(N+2m-2i)\kappa_{\mu^i}/2}) \\
& \cdot \sum_{\nu^{1, \dots, N-1}}
\prod_{k=1}^N (\cW_{\nu^{k-1}\mu^k(\nu^k)^t}(q) q^{\kappa_{\nu^k/2}} Q_k^{|\nu^k|}) \\
& \cdot \sum_{\eta^{1, \dots, N-1}}
\prod_{k=1}^N (\cW_{\eta^{k-1}\mu^{k}(\eta^k)^t}(q) q^{\kappa_{\eta^k/2}} Q_k^{|\eta^k|}).
\end{split}
\end{equation}
This suggests the following partition function:
\begin{equation} \label{def:N>2}
\begin{split}
\tilde{\cZ}^{(m)}_{A_{N-1}} = & \sum_{\mu^1, \dots, \mu^N}
\prod_{i=1}^N (Q_{b_i}^{|\mu^i|} (-1)^{(N+m)|\mu^i|} q^{(N+2m-2i)\kappa_{\mu^i}/2}) \\
& \cdot \sum_{\nu^{1, \dots, N-1}}
\prod_{k=1}^N (W_{\nu^{k-1}\mu^k(\nu^k)^t}(q, e^{-t}) q^{\kappa_{\nu^k/2}} Q_k^{|\nu^k|}) \\
& \cdot \sum_{\eta^{1, \dots, N-1}}
\prod_{k=1}^N (W_{\eta^{k-1}\mu^{k}(\eta^k)^t}(q, e^{t}) q^{\kappa_{\eta^k/2}} Q_k^{|\eta^k|}).
\end{split}
\end{equation}
In other words,
\begin{eqnarray*}
\tilde{\cZ}^{(m)}_{A_{N-1}}
& = & \sum_{\mu^1, \dots, \mu^N}
\prod_{i=1}^N (Q_{b_i}^{|\mu^i|} (-1)^{(N+m)|\mu^i|} q^{(N+2m-2i)\kappa_{\mu^i}/2}) \\
&& \cdot \tilde{\cK}_{\mu^1, \dots, \mu^N}(Q_1, \dots, Q_{N-1}; q, e^{-t})
\cdot \tilde{\cK}_{\mu^1, \dots, \mu^N}(Q_1, \dots, Q_{N-1}; q, e^{t})
\end{eqnarray*}
As in the $N=2$ case one can try to relate it to the generating series of
$$\chi(M(N, k), \Lambda_{-e^{-t}}(\bV_k \otimes \bW_k^*) \otimes
\Lambda_{-e^{-t}}(\bV_k^* \otimes \bW_k) \otimes E_{N, k}^m)(e_1, \dots, e_N, q, q^{-1}).$$

{\bf Acknowledgements}.
This research is partially supported by research grants from NSFC and Tsinghua University.
The final stage of its preparation was carried out in the Center of Mathematical Sciences,
Zhejiang University.
The author thanks Professor Kefeng Liu for his interest in this work.

\end{document}